\documentclass[11pt]{article}
\usepackage{amsmath,latexsym,amssymb,epsfig}
\begin{document}
\setlength{\topmargin}{0in}
\setlength{\headheight}{0.12in}
\setlength{\headsep}{.40in}
\setlength{\parindent}{1pc}
\setlength{\oddsidemargin}{-0.1in}
\setlength{\evensidemargin}{-0.1in}
\headheight 12pt
\headsep 25pt
\footskip 30pt
\textheight  625pt 
\textwidth 170mm
\columnsep 10pt
\columnseprule 0pt
\setlength{\unitlength}{1mm}

\setlength{\parindent}{20pt}
\setlength{\parskip}{2pt}
\newtheorem{thm}{Theorem}[section]
\newtheorem{example}[thm]{Example}
\newtheorem{cor}[thm]{Corollary}
\newtheorem{notation}[thm]{Notation}
\newtheorem{notations}[thm]{Notations}
\newtheorem{definition}[thm]{Definition}
\newtheorem{lemma}[thm]{Lemma}
\newtheorem{claim}[thm]{Claim}
\newtheorem{mthm}[thm]{Meta-Theorem}
\newtheorem{prop}[thm]{Proposition}
\newtheorem{rem}[thm]{Remark}
\newtheorem{conj}[thm]{Conjecture}
\newtheorem{rems}[thm]{Remarks}
\renewcommand {\theequation}{\thesection.\arabic{equation}}
\newcommand{\qed}{\hfill\rule{2mm}{3mm}\vspace{4mm}}
\newcommand{\gdm}{\hfill\vrule  height5pt width5pt \vspace{.1in}}
\def\al{{\alpha}}\def\be{{\beta}}\def\de{{\delta}}
\def\ep{{\epsilon}}\def\ga{{\gamma}}\def\ka{{\kappa}}
\def\la{{\lambda}}\def\om{{\omega}}\def\si{{\sigma}}
\def\ze{{\zeta}}\def\De{{\Delta}}
\def\Ga{{\Gamma}}
\def\La{{\Lambda}}\def\Om{{\Omega}}
\def\th{{\theta}}\def\Th{{\Theta}}
\def\<{\left<}\def\>{\right>}\def\({\left(}\def\){\right)}

\newfam\msbmfam\font\tenmsbm=msbm10\textfont
\msbmfam=\tenmsbm\font\sevenmsbm=msbm7


\scriptfont\msbmfam=\sevenmsbm\def\bb#1{{\fam\msbmfam #1}}

\def\AA{\bb A}\def\BB{\bb B}\def\CC{\bb C}\def\DD{\bb D}
\def\EE{\mathbb{E}}\def\GG{\bb G}\def\HH{\bb H}\def\KK{\bb K}
\def\LL{\bb L}\def\NN{\bb N}\def\PP{\mathbb{P}}\def\QQ{\bb Q}
\def\RR{\mathbb{R}}\def\TT{\bb T}\def\WW{\mathbb{W}}\def\ZZ{\mathbb{Z}}
\def\cA{{\cal A}}\def\cB{{\cal B}}\def\cL{{\cal L}}\def\cC{{\cal C}}
\def\cD{{\cal D}}\def\cF{{\cal F}}
\def\cH{{\cal H}}\def\cI{{\cal I}}\def\cJ{{\cal J}}
\def\cL{{\cal L}}\def\cM{{\cal M}}
\def\cP{{\cal P}}\def\cS{{\cal S}}\def\cT{{\cal T}}\def\cV{{\cal V}}
\def\cW{{\cal W}}\def\cX{{\cal X}}\def\cY{{\cal Y}}\def\cZ{{\cal Z}}
\def\tl{{\tilde \ell}}
\def\tZ{\widetilde{Z}}
\def\hZ{\widehat{Z}}
\def\cz{ {\mathcal U}}
\def\bw{{\mathbf w}}\def\bW{{\mathbf W}}

\def\us{{\underline{s}}}
\def\ut{{\underline{t}}}
\def\uu{{\underline{u}}}

\def\goto{{\rightarrow}}

\title{\large \bf Snake representation of a superprocess\\ in random environment}
\author{Leonid Mytnik\thanks{Faculty of Industrial Engineering,
Technion, Haifa 32000, Israel. 
Partially supported
by a grant from the Israel Science Foundation and by B. and G. Greenberg
Research Fund (Ottawa).}
 \and Jie Xiong\thanks{Department of
Mathematics, University of Tennessee, Knoxville, TN 37996, USA and
Department of Mathematics, Hebei Normal University, Shijiazhuang
050016, China. Research supported partially by NSF DMS-0906907.}
\and Ofer Zeitouni\thanks{Faculty of Mathematics, Weizmann Institute
of Science, Rehovot 76100, Israel and School of Mathematics,
University of Minnesota, MN 55455, USA. Partially supported by NSF
grant DMS-0804133, a grant from the Israel Science
Foundation, and by
the Herman P. Taubman Professorial chair at the Weizmann Institute. } }
\date{January 13, 2011. Revised July 19, 2011.}

\maketitle
\begin{abstract}
    We consider (discrete time) branching particles in a
    random environment which is i.i.d. in time and possibly spatially
    correlated. We prove a representation of the limit process
    by means of a Brownian snake in random environment.
\end{abstract}

\section{Introduction}
\label{sec:1}
\subsection{Superprocesses in random environments}
\label{1.1}
Superprocesses in random environments were
introduced in~\cite{M} as the scaling limits of
particle systems whose branching are
affected by random environments.
In particular the limiting behavior
of the following model has been studied.
At time $t=0$, $K_n\sim n$ particles are located in $\RR^d$. Each
of these $K_n$ particles follows the path of an independent
Brownian motion
until time $t=1/n$. At time $1/n$ each particle independently of the others
either splits into two or dies and then the individual particles in the new population
again  follow the paths of independent Brownian
motions starting at their place of birth,
in the interval $[1/n,2/n)$, and the pattern of
alternating branching and spatial spreading continues.
Let us describe in details the branching mechanism
that was suggested in~\cite{M}.
Let $\{\xi_k(\cdot)\}_{k\geq 0}$ be a sequence of i.i.d.
$\RR^d$-indexed random fields with
mean $0$  and covariance
\[g(x,y)= {\rm Cov}(\xi_k(x),\xi_k(y)),\quad x,y\in \RR^d.\]
At time ${k}/{n}$ each particle,
independently of the others conditionally on $\xi$,
either splits into two with probability
 $$\frac12+\frac{1}{2\sqrt{n}}\xi_k(x)$$
or dies with probability
 $$\frac12-\frac{1}{2\sqrt{n}}\xi_k(x),$$
where $x$ is the location of the
particle. That is, the fields
$\{\xi_k\}_{k\geq 0}$
create the random environment that affects the branching
of the particles.
Define the following measure-valued process
that describes  the evolution of the population:
\begin{eqnarray}
\label{eq:1.1.1}
X^n_t(A)=\frac{{\rm number \; of\; particles \; in \;} A\; {\rm at \; time\;}t}
{n},\; A\subset \RR^d.
\end{eqnarray}
Before proceeding we introduce some notation. 
For a locally compact Polish space $E$, let $\cM_F(E)$
(respectively, $\cM(E)$) be the space of finite (respectively Radon)
non-negative measures on $E$, equipped with the weak  (respectively,
vague) topology~(see Section 3.1 in~\cite{bib:daw91}). In the case
of $E=\RR^d$, we will also write $\cM_F=\cM_F(\RR^d)$ and
$\cM=\cM(\RR^d)$. Both $\mu(\phi)$ and $\langle \phi, \mu\rangle$
denote the integral of a function $\phi$ with respect to measure
$\mu$. For any metric space $E$ let 
$D_E=D_E[0,\infty)$ (resp.
$C_E=C_E[0,\infty)$) be the space of cadlag (resp. continuous)
$E$-valued functions on $[0,\infty)$ endowed with the Skorohod
topology. 
Let
$\cC^k(\RR^d)$ (resp. $\cC^k_b(\RR^d)$) be the set of continuous
(resp. bounded continous) functions with continuous (resp. bounded
continuous) partial derivatives of order $k$ or less. Also we define
$\cB(\RR^d)$ to be the set of bounded measurable functions on
$\RR^d$.

It was
shown in~\cite{M}, under
some additional technical assumptions on $\xi$, that if
$$ X^n_0 \Rightarrow X_0=:\mu\,, \;\; {\rm in}\; \cM_F\,,$$
then
$$ X^n \Rightarrow X, \;\; {\rm in}\; D_{\cM_F}[0,\infty)\,.$$
Here $X$ is a process in $C_{\cM_F}[0,\infty)$ which is the unique solution to the following martingale problem:
 $\forall\
\phi\in \cC^2_b(\RR^d)$,
\begin{equation}\label{MP-01}
M^{\phi}_t\equiv \left<X_t,\phi\right>-\left<\mu,\phi\right>-
\frac12 \int^t_0\left<X_s,\De\phi\right>ds,\qquad t\ge 0
\end{equation}
is a continuous martingale with quadratic variation process
\begin{eqnarray}\label{MP-02}
\left<M^{\phi}\right>_t&=&\int^t_0\left<X_s,\phi^2\right>ds\\
&&+\int^t_0\int_{\RR^d}\int_{\RR^d}g(x,y)\phi(x)\phi(y)X_s(dx)X_s(dy)ds,\qquad
t\ge 0.\nonumber
\end{eqnarray}

In this paper we introduce some minor changes into the above model.
Instead of the binary branching  we  assume that each particle gives
birth to a number of particles distributed according to the
geometric distribution with parameter
$\frac12-\frac{1}{4\sqrt{n}}\xi_k(x)$; that is, if $N$ is the number
of offspring of the particle located at $x$ at time $k/n$, then
\begin{eqnarray}
\label{eq:1.2} \PP(N=m|\xi)=
\left(\frac12+\frac{1}{4\sqrt{n}}\xi_k(x)\right)^m\left(\frac12-
\frac{1}{4\sqrt{n}}\xi_k(x)\right),\;\;\;m=0,1,2,\ldots.
\end{eqnarray}
In particular,
conditioned on the environment
$\xi$, the expected number of offspring of a particle at $x$ at time $k/n$
is
\begin{eqnarray}
\label{eqm32}
\frac{\frac12+\frac{1}{4\sqrt{n}}\xi_k(x)}{\frac12-
\frac{1}{4\sqrt{n}}\xi_k(x)}
=1+\frac{1}{\sqrt{n}}\xi_k(x)+\frac{1}{2n}\xi_k(x)^2+o\(\frac1n\).
\end{eqnarray}

Compared with \cite{M},
we also allow $\xi$ to be  slightly
more general, that is,  we assume that $\{\xi_k(\cdot)\}_{k\geq 1}=\{\xi_k^n(\cdot)\}_{k\geq 1}$ is
a sequence of i.i.d. random fields with
mean $\nu/\sqrt n$, for some $\nu\in \RR$,  and covariance
\begin{equation}\label{eq-Nov26a}
    g(x,y)= {\rm Cov}(\xi_k(x),\xi_k(y)),\quad x,y\in \RR^d.
\end{equation}

Let $X^n$ be defined for this model as in (\ref{eq:1.1.1}).
  By the same
argument as in~\cite{M} one can prove that the limit of $\{X^n\}_{n\geq 1}$
is the (unique) solution to the following martingale problem: $\forall\
\phi\in {\cal C}^2_b(\RR^d)$,
\begin{equation}\label{MP-1}
M^{\phi}_t\equiv \left<X_t,\phi\right>-\left<\mu,\phi\right>-
\frac12 \int^t_0\(\left<X_s,\De\phi\right>+
\<X_s,(\nu+\bar g/2)\phi\>\)ds,\qquad
t\ge 0
\end{equation}
is a continuous martingale with quadratic variation process
\begin{eqnarray}\label{MP-2}
\left<M^{\phi}\right>_t&=&2\int^t_0\left<X_s,\phi^2\right>ds\\
&&+\int^t_0\int_{\RR^d}\int_{\RR^d}g(x,y)\phi(x)\phi(y)X_s(dx)X_s(dy)ds,\qquad
t\ge 0\nonumber
\end{eqnarray}
where $\bar{g}(x)=g(x,x)$.
(Note the different factor multiplying the term $\<X_s,\phi^2\>$ here compared 
with
\eqref{MP-02}, which comes from the different variances of the geometric and Bernoulli distributions.)

\subsection{Brownian snake}
\label{1.2}
The main purpose of this
 paper is to study the Brownian snake representation of the process that solves
the above martingale problem~(\ref{MP-1}-\ref{MP-2}). For a nice
introduction into the topic the reader is referred to~\cite{LeG99}.
The classical Brownian snake was used to study different
properties of super-Brownian motion. Loosely speaking if
$\{\WW_s\}_{s\geq 0}$ is a Brownian snake then for each $s\geq 0$,
 $\WW_s$ is a stopped Brownian path.
 To be more precise we call the pair
 $w=(\bw,\zeta)\in
C_{\RR^d}[0,\infty)\times
 \RR_+$ a stopped path in $\RR^d$ if
 for each $t\geq \zeta$, $\bw(t)=\bw(\zeta)$.
 $\zeta$ is called the lifetime of the
 path $w$ and sometimes is denoted by $\zeta_w$ or $\zeta(w)$.
 Let $\cW$ denote the space of all
 stopped paths in $\RR^d$ equipped with the distance
\[ d(w,w')= \sup_{t\geq 0}|\bw(t)-\bw'(t)|+ |\zeta_w-\zeta_{w'}|. \]
We will also use the notation 
$\hat
w=\bw(\zeta_w)$ for the terminal point of $w$. For any $x\in \RR^d$
we denote by $\bar x$ the path with lifetime $0$ constantly equal to
$x$. 
If $w=(\bw,\zeta)$ is a stopped path
then with some abuse of notation we will sometimes set
$w(s)=\bw(s)$ for any $s\geq 0$.

The usual  Brownian snake can be thought of
as a limit of the so-called discrete snakes that we will now define.
Let $\{Y_{k/n^2}^n\}_{k=0,1,\ldots}$ be a rescaled simple random walk on $\ZZ_+/n$ reflected at the origin, that is, the time between the steps is
$1/n^2$ and the size of the jump is $\pm 1/n$ with equal probabilities. Explicitely,
\begin{eqnarray*}
\PP\(Y^n_{(k+1)/n^2}-Y^n_{k/n^2}=\pm
1/n\)&=&\frac12,\;\;{\rm if}\; Y^n_{k/n^2}\geq 1/n,\;  k=0,1,\ldots,\\
Y^n_{(k+1)/n^2}&=&1/n,\;\;{\rm if}\; Y^n_{k/n^2}= 0,\;  k=0,1,\ldots.
\end{eqnarray*}
We also let $Y^n_\cdot$ 
be
constant between the jumps. The process $Y^n_\cdot$ is called the
contour or lifetime process of the discrete snake $\WW^n$. That is
for each $s\geq 0$, the snake $\WW^n_s=(\bW^n_s, Y^n_s)$
is a stopped
path with life time $Y^n_s$. We next define the paths of the  snake.
Fix $x\in \RR^d$ and set
$$ \WW^n_0=\bar x.$$
Let $\eta_1,\eta_2,\ldots$ be a sequence of independent Brownian
paths stopped at time $1/n$, independent of the contour process. Let
$\WW^n_{k/n^2}=(\bW^n_{k/n^2}, Y^n_{k/n^2})$ 
be the stopped path at time
$k/n^2$ with lifetime
 $\zeta_{\WW^n_{k/n^2}}=Y^n_{k/n^2}$. Then define
\begin{eqnarray}
\label{eq:1.8}
\bW^n_{(k+1)/n^2}(\cdot)=\left\{ \begin{array}{rcl}
    \bW^n_{k/n^2}(\cdot \wedge (Y^n_{k/n^2}-1/n)),&&{\rm if}\;
Y^n_{(k+1)/n^2}=Y^n_{k/n^2}-1/n,\\
  \bW^n_{k/n^2}\odot \eta_k(\cdot),&& {\rm if}\; Y^n_{(k+1)/n^2}=
Y^n_{k/n^2}+1/n,
\end{array}
\right.
\end{eqnarray}
where $\eta_1\odot\eta_2$ denotes  the
concatenation of two paths $\eta_1$ and $\eta_2$ in the obvious way.
In words, if the lifetime $Y^n$ goes down by $1/n$
we erase the path of the snake from the tip by $1/n$, or to put it
differently, we reduce its lifetime by $1/n$.
If $Y^n$ goes up by $1/n$ we add the path $\eta$ to the tip of the
snake. Then 
we define
\[ \WW^n_{(k+1)/n^2}(\cdot)= (\bW^n_{(k+1)/n^2}(\cdot), Y^n_{(k+1)/n^2}(\cdot)).\]
 This way we constructed a sequence of discrete snakes. As is the case
 for $Y^n_\cdot$,
 we define
$\WW_s^n(\cdot)=\WW_{\lfloor sn^2\rfloor/n}^n(\cdot)$.
The sequence of processes $\WW^n$  converges, as $n\rightarrow\infty$,
to a continuous time
Brownian snake (see e.g. Proposition 2.2 in~\cite{LeG96}).

We next
describe  the connection between the snake process and the branching
Brownian motion.
Define the discrete version of the local time as the
rescaled number of upcrossings
  of $Y^n$ from the corresponding level:
\begin{eqnarray}
\label{eq:1.3}
\ell^{n,m/n}_{s}=n^{-1}\sum^{\lfloor sn^2\rfloor}_{i=0}
{\bf 1}_{\{Y^n_{i/n^2}=m/n,Y^n_{(i+1)/n^2}=(m+1)/n\}}.
\end{eqnarray}
We 
also define, for
$t\geq 0$,
\begin{eqnarray}
\label{eq:30_12_10}
\ell^{n,t}_{s}=\ell^{n,\lfloor tn\rfloor/n}_{s}\,.
\end{eqnarray}
Since 
$s\mapsto \ell^{n,t}_{s}$ is increasing we define the measure
$\ell^{n,t}(ds)$ in an obvious way. In 
fact this convention
will be used throughout the paper: for any non-decreasing function
$r\mapsto f_r$ on $\RR_+$, $f(dr)$, with a slight abuse of notation,
will denote the corresponding measure defined via
$f((a,b])=f_b-f_a$, for any $b\geq a$.

For any $a\geq 0$ 
introduce the inverse local time at level $a$ as 
\begin{eqnarray}
\label{eq:1.4}
\tau^{n,a}_r=\frac{1}{n^2}\inf\{k: \ell^{n,a}_{k/n^2}> r\}.
\end{eqnarray}
\begin{figure}[t]
\epsfig{file=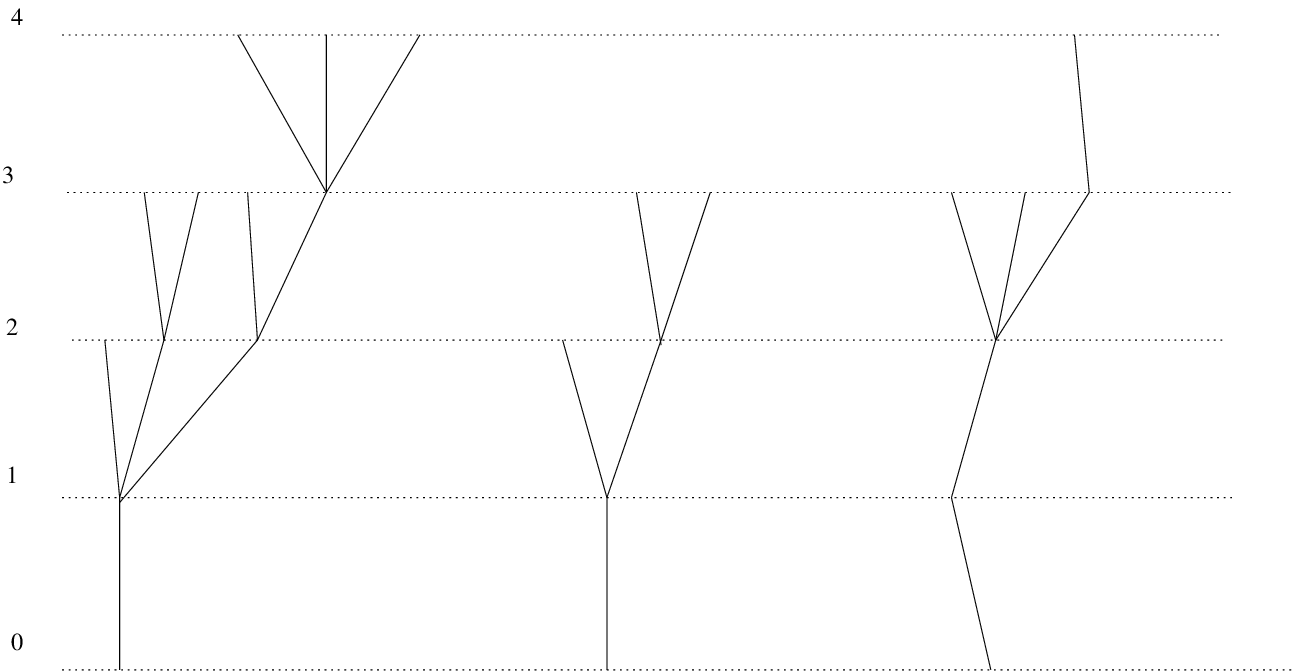,height=6.5cm, width=6.5cm, angle = 0}
\epsfig{file=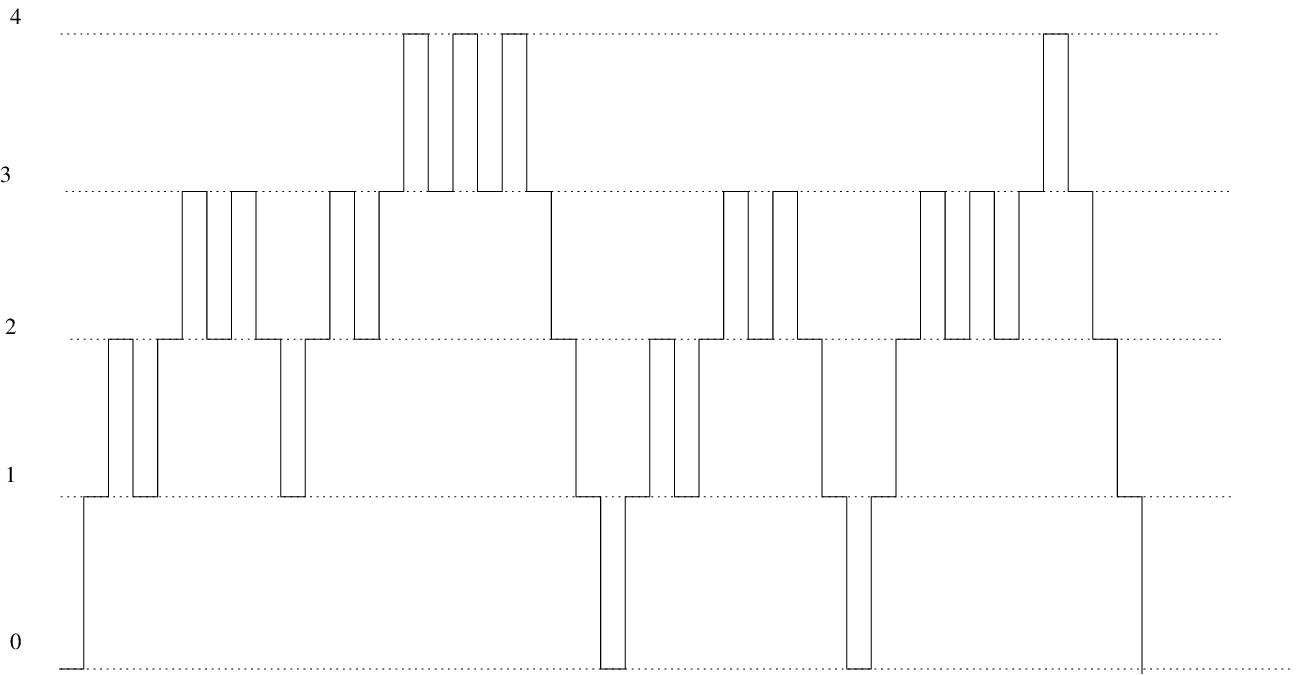,height=6.5cm, width=8.5cm, angle = 0}
\caption{Geneology of the particle system (left) and contour
process (right). In this picture, $\tau^{1,0}_1=23 $, $\tau^{1,0}_2=33$ and
$\tau^{1,1}_4=26$.}
\end{figure}

For any $a\geq0$ and $r_1<r_2$  define the measure valued process
$X^{n,r_1,r_2}_{a,t}$ so that, for $\phi\in \cB(\RR^d)$
\begin{eqnarray}
\label{eqm40}
X^{n,r_1,r_2}_{a,t}(\phi)&\equiv&\int_{\tau^{n,a}_{r_1}}^{\tau^{n,a}_{r_2}}
\phi(\WW^n_s(Y^n_s))\ell^{n,t}(ds)\,,\; t\geq a.
\end{eqnarray}
It is easy to see that
$X^{n,r_1,r_2}_{a,k/n},\;\;k\geq \lfloor an\rfloor$,
is the measure-valued process constructed
 in the previous section starting at ``time''
 $a_n=\lfloor an\rfloor/n$
 such that
 \[ X^{n,r_1,r_2}_{a,a_n}(1)=X^{n,r_1,r_2}_{a_n,a_n}(1)=r_2-r_1\,,\]
and therefore
\begin{eqnarray}
X^{n,r_1,r_2}_{a,\cdot}\Rightarrow X^{r_1,r_2}_{a,\cdot},\;\;
\end{eqnarray}
where $X^{r_1,r_2}_{a,\cdot}$
solves the martingale problem  starting at time $a$
such that
\[X^{r_1,r_2}_{a,a}(1)=r_2-r_1\]
and,
$\forall\
\phi\in {\cal C}^2_b(\RR^d)$,
\begin{equation}\label{MP-3}
    M^{\phi}_{a,t}
    \equiv \left<X_{a,t}^{r_1,r_2},\phi\right>-
    \left<X_{a,a}^{r_1,r_2},\phi\right>-
\frac12\int^t_a\left<X_{a,s}^{r_1,r_2},\De\phi\right>ds,\qquad
t\ge a,
\end{equation}
is a continuous martingale with quadratic variation process
\begin{eqnarray}\label{MP-4}
    \left<M^{\phi}_a\right>_t&=&2\int^t_a
    \left<X_{a,s}^{r_1,r_2},\phi^2\right>ds, \qquad
t\ge a.
\end{eqnarray}

\subsection{Our model}
We finally
define the discrete snake in random environment corresponding to
the branching processes in random environment described in Section~\ref{1.2}.
The
main difference with the ``fixed environment''
case is that here the snake cannot be
constructed conditionally on the lifetime process.
Both processes have to be constructed
simultaneously.

The environment $\{\xi_k(\cdot)\}_{k\geq
0}=\{\xi_k^n(\cdot)\}_{k\geq 0}$ is assumed to consist of a sequence
of i.i.d. random fields, satisfying  $|\xi^n_k(x)|\leq \sqrt{n}/2$
and 
$\sup_n
\EE(|\xi_k^n(x)|^3)<\infty$, with mean $\nu/\sqrt n$, for some
$\nu\in \RR$, and covariance $g(x,y)$ as in (\ref{eq-Nov26a}), with
$\|\bar g\|_\infty<\infty$.

Now define the snake
with lifetime processes $\WW^n=(\bW^n, Y^n)$ as follows. Fix a
constant $K_1>0$. Let $Y^n_0=0$ and $\WW^n_0=\bar{x}$ with $x\in
\RR^d$. Suppose we are given $(\bW^n_{k/n^2}, Y^n_{k/n^2})$ for some
$k\geq 0$.
 $(\bW^n_{(k+1)/n^2}, Y^n_{(k+1)/n^2})$ will be defined as follows. If
 $Y^n_{k/n^2}\not\in \{0,K_1\}$, then
conditionally on $\xi$ and  $(\bW^n_{l/n^2}, Y^n_{l/n^2}), l\leq k$ we set
\[\PP\(Y^n_{(k+1)/n^2}-Y^n_{k/n^2}=\pm
1/n |\xi,  \bW^n_{l/n^2}, Y^n_{l/n^2}, l\leq k\)=
 \frac12\pm\frac{1}{4\sqrt{n}}\xi_{Y^n_{k/n^2}}\(\hat{\WW}^n_{n^{-2}k}\),\]
where
we introduced 
above the  notation for the ``tip'' of the snake:
\[ \hat{\WW}^n_{k/n^2} =  \bW^n_{k/n^2}\(Y^n_{k/n^2}\).\]
If $Y^n_{k/n^2}= 0$, then with probability one we set $Y^n_{(k+1)/n^2}=1/n$.
If
$Y^n_{k/n^2}=K_1$, then with probability one we set
$Y^n_{(k+1)/n^2}=K_1-1/n$.
(That is,
the process is reflected at height $K_1$; a
similar approach
of introducing a super-critical branching mechanism
via a
reflection of the lifetime process was used by J.-F. Delmas in~\cite{D08}.)


Let $\eta_1,\eta_2,\ldots$ be  a sequence of independent
Brownian motions stopped at time $1/n$.
Given the evolution of the lifetime
process $Y^n$ until time $(k+1)/n^2$, the path of the
Brownian snake $\WW$ at time $(k+1)/n^2$ is defined exactly as in
(\ref{eq:1.8}).

We next explain the connection between the snake and branching
particle system in random environment which is analogous to the
connection that exists between the processes in a constant
environment. Define the rescaled local time $\ell^{n,t}_s$ for $Y^n$
as in (\ref{eq:1.3}), (\ref{eq:30_12_10}) and the inverse local time
as in (\ref{eq:1.4}). For 
any $r_1,r_2>0, a\geq 0$, we define the
measure-valued process in the same way as it is done
in~(\ref{eqm40}):
\begin{eqnarray}
\label{eq:31_12_01}
X^{n,r_1,r_2}_{a,t}(\phi)&\equiv&\int_{\tau^{n,a}_{r_1}}^{\tau^{n,a}_{r_2}}
 \phi(\bW^n_s(Y^n_s))
\ell^{n,t}(ds),\;\;\; t\geq a,
\end{eqnarray}
for all $\phi\in \cB(\RR^d)$.
This process characterizes the branching particle picture in random environment
with offspring distribution given by~(\ref{eq:1.2}) and
starting  with
$\lfloor (r_2-r_1)n\rfloor$ particles at the site $x\in \RR^d$ at time $t=a$. In the case of $r_1=0, r_2=r, a=0$, we will use the notation
\begin{eqnarray}
\label{eq:Oc5_1}
X^{n,r}_{0,t}&\equiv&X^{n,r_1,r_2}_{0,t},\;\;\; t\geq 0,
\end{eqnarray}
for the corresponding process.

The following is our first main result.
\begin{thm}
\label{thm:1} Fix $K_1>0$. Then the sequence of processes
$\{\WW^n\}_{n\geq 1}$ $=$ $\{(\bW^n,Y^n)\}_{n\geq 1}$ is $C$-tight
in $D_{\cW}$. 
Let
$\WW=(\bW,Y)$ be an arbitrary limiting point, let $\ell^a$ be a
local time of $Y$ at  level $a$ and let $\tau^a(r)$ be the inverse
of the local time.
 Fix an arbitrary $r>0$.  Then
\begin{eqnarray}
\label{20_2}
X_t^r(\phi) = \int_0^{\tau^0_r}
\phi(\hat{\WW}_s) \ell^t(ds), \;\; \phi\in \cB(\RR^d), \;
t\in [0,K_1],
\end{eqnarray}
is the  measure-valued process satisfying the
martingale problem (\ref{MP-1}-\ref{MP-2}) on $[0,K_1]$,
with
 $X_0^r=r\delta_x$.
\end{thm}
(The uniqueness of solutions to \eqref{MP-1}--\eqref{MP-2} can be read off 
\cite{M}, however the uniqueness of the snake $\WW$ is an open problem.)

In the particular case of a spatially ``smooth" 
random environment we can give another description of
the snake process. It is easy to check from our assumptions on
$\xi^n$ that if we define
\begin{eqnarray}
B^n_s(x)\equiv \frac{1}{\sqrt n}\sum_{i=1}^{\lfloor sn\rfloor} \xi^n_i(x),
\end{eqnarray}
then
\[ B^n \rightarrow B,\]
where
$\frac{\partial B_t(y)}{\partial t}$ is a Gaussian generalized noise on $\RR_+\times\RR^d$, white in time
and colored in space, such that
\begin{eqnarray}
\label{whitefield}
   \EE( B_t(x))&=&t\nu,\;\;\forall t\geq 0,\nonumber\\
{\rm Cov}\left(\frac{\partial B_t(x)}{\partial t},\frac{\partial B_s(y)}{\partial s}\right) &=& \delta_0(t-s)g(x,y),\\
 B_0&=&0,
\nonumber
\end{eqnarray}
where 
$\delta_0(\cdot)$ is the Dirac measure at $0$. Given the
result on the tightness of $\{\WW^n\}_{n\geq 1}$, one can easily
deduce that the pair $\{(\WW^n,B^n)\}_{n\geq 1}$ is tight. In what
follows we assume that $(\WW,B)$ is a limit point of the tight
sequence $\{(\WW^n,B^n)\}_{n\geq 1}$, and we recall that
$\WW=(\bW,Y)$.

Our aim is to introduce a particular
functional of the limiting snake that has
a simple
semimartingale decomposition.
The definition of the functional is motivated by the one used
by Dhersin and Serlet \cite{DS} and also by a functional used to transform
Brox's diffusion
into a martingale, see \cite{Shi}.
For $w\in\cW$, let
\[F(w)=\int^{\ze(w)}_0 e^{-B_r(w(r))} dr.\]

Our second main result is the following.
For technical reasons, we restrict attention to branching laws constructed
directly from smooth Gaussian fields.
\begin{thm}
\label{thr:5.1} Fix $K_1>0$. Let
$\tilde B$ be a Gaussian field  as in \eqref{whitefield} so that 
$\tilde B\in
C_{\cC^{2}(\RR^d)}[0,\infty)$, 
a.s.., and set
\begin{equation}
\label{eq-newwhitefield}
\frac{\xi_j^n(y)}{\sqrt{n}}=
\left(\tilde B_{\frac{j}{n}}(y)-\tilde B_{\frac{j-1}{n}}(y)\right)
{\bf 1}_{\{
|\tilde B_{\frac{j}{n}}(y)-\tilde B_{\frac{j-1}{n}}(y)|
<1/2\}}\,.
\end{equation}
Let $(\WW,B)$ denote a limit point of $(\WW^n,B^n)$.
Then $B$ has the distribution of $\tilde B$, and
there
exists a Brownian motion $\beta$ such that
\begin{eqnarray}
\nonumber
F(\WW_t)&=& \int_0^t e^{-B_{Y_s}(\hat{\WW}_s)}\left\{-\frac{1}{2}\Delta B_{Y_s}(\hat{\WW}_s)+\frac{1}{2}\sum_{i=1}^d
  \left(\frac{\partial}{\partial x_i} B_{Y_s}(\hat{\WW}_s)\right)^2\right\}\,ds \\
\label{eq:5.1}
&&+ \ell^0_t -\int_0^t e^{-B_{K_1}(\hat{\WW}_s)}\,\ell^{K_1}(ds)+
 \int_0^t e^{-B_{Y_s}(\hat{\WW}_s)}\,d\beta_s\,.
\end{eqnarray}
\end{thm}
\begin{rem}
The first term on the right side of~(\ref{eq:5.1}) can be written as
$$ \int_0^t \frac{1}{2}\Delta_x e^{-B_{Y_s}(x)}|_{x=\hat{\WW}_s}\,ds, $$
and it comes from the fact that $\bW_s(\cdot)$ 
is a Brownian path.
\end{rem}
\begin{rem}
It is plausible that one may relax the assumptions on the fields 
$\{\xi^n\}$ in Theorem 
\ref{thr:5.1}; however, some strong assumptions that ensure good spatial
approximation of $B$ by partial sums of the $\xi^n$s seem to be 
crucial.
\end{rem}
Note that in the case of constant function $g$, for every $s$, $B_s(\cdot)$ is
a constant
function in space, and hence we immediately have the following corollary, which
for simplicity we state only in case $\nu=0$.
A similar result (without the reflection) can be found in
\cite{Sei}.
\begin{cor}
\label{cor-brox}
 Let $g\equiv1$ and $\nu=0$. Then $Y$ is the Brox diffusion
reflected at $0$ and $K_1$.
\end{cor}
See the appendix for the definition of the Brox diffusion.

\subsection{Structure of the paper}
In the next section, we derive some standard estimates on survival probability
for branching processes in a random environment. Section
\ref{sec-tcont} is concerned with the proof of tightness of the contour process.
(Because of dependence through the environment, natural arguments involving
stopping times such as Aldous' tightness criterion cannot be applied directly,
and extra care has to be employed in separating dependence on the level
of the contour process from dependence on the lifetime of the process.)
Most of the work is devoted to proving that large upward jumps of the contour
process are unlikely; downward jumps are then handled by a time
reversal argument.
Section  \ref{sec-4} is devoted to the proof of tightness of the snake process
and its local time, and a completion of the proof of
Theorem \ref{thm:1}.
Section \ref{sec-5} is devoted to the description of the snake provided
in Theorem \ref{thr:5.1}, while the appendix is devoted to the description
of the contour process for
environments with no spatial dependence,
providing in particular a direct proof of
Corollary \ref{cor-brox}, that bypasses the need to consider the Brownian snake.

\noindent
{\bf Notation} Throughout, $C,K$ denote  generic constants whose values
may change from line to line. Numbered constants (such as $K_1,c_0,C_m,\delta_{4.2}$, etc.) are fixed and do not change throughout the paper.

\section{Asymptotics for survival probability and useful bounds}
\label{sec-asymp} \setcounter{equation}{0} We start with a lemma
that describes the asymptotics for survival probability for
classical branching processes. For any $n\geq 1$ let $\{M^n_l\,,
l=0,1,2,\ldots\}$ be the branching process with geometric offspring
distribution with parameter
\[ p=1/2 -b_n/4n\]
for some $b_n\in (-2n,2n)$.
 That is if $Z^n$ is the number of offspring
in the process $M^n$, then
\[ \PP(Z^n=k) = p(1-p)^k,\;\; k=0,1,2,\ldots.\]

For $\delta>0$ define
\[ h(b,\delta)  = \left\{
\begin{array}{l}
\frac{b}{1-e^{-b\delta}},\;\; b\not= 0,\\
\frac{1}{\delta},\;\;b=0.
\end{array}
\right.
\]
\begin{lemma}
\label{lem:3}
Assume
\[ \lim_{n\rightarrow \infty} b_n=b,\]
and $M^n_0=1$ for all $n\geq 1$.
Then for any $\delta>0$,
\[ \lim_{n\to \infty} n \PP(M^n_{[n\delta]}>0)= h(b,\delta).
\]
\end{lemma}
\paragraph{Proof:}
For $b=0$ the result  is well-known (see e.g. \cite[Theorem II.1.1]{Perkins}
for a more general result).
 While we believe that the result is also known for
$b>0$, we were unable to locate a reference
and thus for the sake of completeness we provide a proof.

Let $f(s)$ be the generating function of $Z^n$, that is
\[ f(s)=\frac{1/2-b_n/4n}{1-(  1/2+b_n/4n)s},\;\; 0\leq s\leq 1.\]
Define $f_0(s)=s, \; f_1(s) = f(f_0(s))=f(s)$ and in general
\[ f_k(s)= f(f_{k-1}(s)),\;\; 0\leq s.\]
Then by the branching property,
\[ E[ s^{M^n_k}| M^n_0=1] =f_k(s).\]
Therefore,
\begin{eqnarray*}
\PP(M^n_k=0| M^n_0=1) &=& f_k(0)
=f(f_{k-1}(0))\\
&=& \frac{1/2- b_n/4n}{1-(1/2+b_n/4n)f_{k-1}(0)}
.\end{eqnarray*}
Fix $k=\lfloor n\delta\rfloor$ and define
\[ y_l=  k (1-f_l(0)),\;\; l=0,1,\ldots, k. \]
One has
\begin{eqnarray}
y_l &=& \frac{ (1/2+b_n/4n)y_{l-1}}{1-(1/2+b_n/4n)(1-y_{l-1}/k)}
.\end{eqnarray}
 Let
\[ z_{l}=1/y_l\,, l=0,\ldots,k\,.\]
Then
$z_0=1/k$
and
\[ z_l = z_{l-1}d_n +1/k,\;\; l=1,\ldots,k,\]
with
\begin{eqnarray}
d_n&=&  \frac{1/2-b_n/4n}{1/2+b_n/4n}
.\end{eqnarray}
By iterating we get
\begin{eqnarray}
z_k&=& \frac{1}{k}\left(d_n^k + \frac{1-d_n^k}{1-d_n}\right)
.\end{eqnarray}
Now as $n\rightarrow \infty$ we have
\[ 1-d_n = \frac{b_n}{n(1+\frac{b_n}{2n})} \sim \frac{b_n}{n}.\]
Also recall that $k=\lfloor n\delta\rfloor$ and hence
$$d_n^k=   \left(\frac{1/2-b_n/4n}{1/2+b_n/4n}\right)^{\lfloor n\delta\rfloor}
\sim (1-b_n/n)^{ n\delta}
\sim e^{-b_n\delta}.
$$
Therefore
$$\lim_{n\rightarrow \infty} z_{\lfloor n\delta\rfloor}=
 \lim_{n\rightarrow \infty}
 \frac{e^{-b_n\delta}+\frac{1-  e^{-b_n\delta}}{b_n/n}}{n\delta}
= \frac{1-  e^{-b\delta}}{b\delta}
.$$
Since $n \PP(M_{\lfloor n\delta\rfloor]}>0)\sim
 \frac{y_{\lfloor n\delta\rfloor}}{\delta}=
\frac{1}{\delta z_{\lfloor n\delta\rfloor}}$, this concludes the proof.
\gdm

Returning to the random environment case,
let $\tilde{M}^n$  denote the total mass  of the branching Brownian
motion in random environment  $X^n$ (with geometric offspring
distribution) defined in Section~\ref{1.1} (that is, $
n^{-1}\tilde{M}^n_k = \langle X^n_{k/n}, 1\rangle$) and let $M^n$
be as above with
\begin{eqnarray}
\label{equt:nov6}
b=\nu + \Vert\bar g\Vert_{\infty}/2.
\end{eqnarray}
\begin{rem}
    \label{rem-new}
    With our assumptions, it is easy to see
that the $m$-th
moment (for any $m\geq 2$) of the absolute value of the expected
(conditioned in the environment)
number of offspring minus 1
of a particle at $x$
at time $k/n$, see (\ref{eqm32}), is bounded by
$C_m/n$,  for an appropriate constant $C_m$. Moreover, the absolute
value of the first moment of
the number of offspring minus 1
of a particle at $x$
at time $i/n$, see (\ref{eqm32}), is bounded by
$C_1/n$,  for an appropriate constant $C_1$.
\end{rem}

\begin{lemma}
\label{lem:4}
Let $\tilde{M}^n_0=1$ for all $n\geq 1$. Then
\[ \limsup_{n\rightarrow \infty} n \PP(\tilde{M}^n_{\lfloor n\delta\rfloor} >0)\leq  h(b,\delta).\]
\end{lemma}
\paragraph{Proof:}
For $i=1,2,\ldots, \tilde{M}^n_k$ we denote by
$\cz_{i,k}(t), t\in [k/n, (k+1)/n],$  the position at time $t$
of the $i$-th particle that was born at time $k/n$.
That is we have
\[ X^n_{k/n} =\frac{1}{n}\sum_{i=1}^{\tilde{M}^n_k} \delta_{\cz_{i,k}(k/n)}\,.\]
Moreover if $Z^n_{i,k+1}$ is the number of offspring at
time $(k+1)/n$ of the $i$-th
particle that was born at time $k/n$, then we also have
\[ X^n_{k/n} =\frac{1}{n}\sum_{i=1}^{\tilde{M}^n_{k-1}} Z^n_{i,k}
\delta_{\cz_{i,k-1}(k/n)}\,.\]
We write for simplicity
$ \xi_{i,k}=\xi_k(\cz_{i,k-1}(k/n))$ and denote by $\cF_k^\xi$
the sigma-algebra generated by the environment $\{\xi_j(\cdot),j\leq k\}$.
We have, 
for $s\in (0,1)$,
\begin{eqnarray*}
\EE\left(s^{\tilde{M}^n_k}\right)&=& \EE\left( \Pi_{i=1}^{\tilde{M}^n_{k-1}} s^{Z^n_{i,k}}\right)\\
\nonumber &=&\EE\left( \Pi_{i=1}^{\tilde{M}^n_{k-1}} \EE\left(
\left. s^{Z^n_{i,k}}\right| X^n_{\frac{k}{n}-}, {\cal F}^{\xi}_{k}\right)\right)
\\
\nonumber
&=& \EE\left( \left(\frac{1}{2-s}\right)^{\tilde{M}^n_{k-1}} \Pi_{i=1}^{\tilde{M}^n_{k-1}}
\frac{1-\frac{\xi_{i,k}}{2\sqrt n}}{1-\frac{\xi_{i,k}s}{2(2-s)\sqrt n}}
 \right)
\,.\end{eqnarray*}
Therefore,
\begin{eqnarray}
\EE\left(s^{\tilde{M}^n_k}\right)&=& 
 \EE\left( \left(\frac{1}{2-s}\right)^{\tilde{M}^n_{k-1}}
\EE\left(\left. \Pi_{i=1}^{\tilde{M}^n_{k-1}}
\frac{1-\frac{\xi_{i,k}}{2\sqrt
n}}{1-\frac{\xi_{i,k}s}{2(2-s)\sqrt n}}\right| X^n_{\frac{k}{n}-}\right)
 \right)\\
\label{equt:03_3}
&\geq&  \EE\left( \left(\frac{1}{2-s}\right)^{\tilde{M}^n_{k-1}}
 \exp\left({
 \EE\left(\left.\sum_{i=1}^{\tilde{M}^n_{k-1}}
\log\left(\frac{1-\frac{\xi_{i,k}}{2\sqrt
n}}{1-\frac{\xi_{i,k}s}{2(2-s)\sqrt n}}\right) \right|
X^n_{\frac{k}{n}-}\right]}\right)
\right),
\nonumber
\end{eqnarray}
where the last inequality follows by Jensen inequality.  Since
$|\xi_{i,k}|<\sqrt{n}/2$
 we get  by trivial estimates that for $n$ large enough
\begin{eqnarray}
\nonumber
 \frac{1-\frac{\xi_{i,k}}{2\sqrt n}}{1-\frac{\xi_{i,k}s}{2(2-s)\sqrt n}}
&\geq&\left( 1-\frac{\xi_{i,k}}{2\sqrt n}\right)\left(1+\frac{\xi_{i,k}s}{2(2-s)\sqrt n}+
  \left(\frac{\xi_{i,k}s}{2(2-s)\sqrt n}\right)^2 +  \left(\frac{\xi_{i,k}s}{2(2-s)\sqrt n}\right)^3 \right)
\\
\nonumber
&=&1+\frac{\xi_{i,k}}{2\sqrt n}\left(\frac{s}{2-s}-1\right)+
   \frac{\xi_{i,k}^2s}{4(2-s) n}\left(\frac{s}{2-s} -1\right)\\
\nonumber
&&\mbox{}+   \left(\frac{\xi_{i,k}s}{2(2-s)\sqrt n}\right)^2\frac{\xi_{i,k}}{2\sqrt n}\left(\frac{s}{2-s} -1\right)
\\
\label{equt:03_1}
&\geq&1-\frac{\xi_{i,k}(1-s)}{(2-s)\sqrt n}-
\frac{\xi_{i,k}^2s(1-s)}{2(2-s)^2 n}-
c_{\ref{equt:03_1}}(1-s)n^{-3/2}|\xi_{i,k}|^3.
\end{eqnarray}
Again by trivial estimate on the logarithmic function we get
\begin{eqnarray}
\nonumber
\log\left(\frac{1-\frac{\xi_{i,k}}{2\sqrt n}}{1-\frac{\xi_{i,k}s}{2(2-s)\sqrt n}}\right)&\geq&
\log\left(1-\frac{\xi_{i,k}(1-s)}{(2-s)\sqrt n}-
   \frac{\xi_{i,k}^2s(1-s)}{2(2-s)^2 n}-
   c_{\ref{equt:03_1}}(1-s)n^{-3/2}|\xi_{i,k}|^3\right)\\
&\geq&-\frac{\xi_{i,k}(1-s)}{(2-s)\sqrt n}-
   \frac{\xi_{i,k}^2s(1-s)}{2(2-s)^2 n}
\nonumber
\\
&&\mbox{}
-\frac{\xi^2_{i,k-1}(1-s)^2}{2(2-s)^2n}-
c_{\ref{equt:03_2}}(1-s)n^{-3/2}
|\xi_{i,k}|^3
\nonumber
\\
\label{equt:03_2}
&=&-\frac{\xi_{i,k}(1-s)}{(2-s)\sqrt n}-
   \frac{\xi_{i,k}^2(1-s)}{2(2-s)^2 n}
-c_{\ref{equt:03_2}}(1-s)n^{-3/2}
|\xi_{i,k}|^3
,\end{eqnarray}
for all $n$ sufficiently large.
Take an expectation to get
\begin{eqnarray}
\lefteqn{\EE\left( \left. -\frac{\xi_{i,k}(1-s)}{(2-s)\sqrt n}-
   \frac{\xi_{i,k}^2(1-s)}{2(2-s)^2 n}
-c_{\ref{equt:03_2}}(1-s)n^{-3/2 } |\xi_{i,k}|^3\right| X^n_{\frac{k}{n}-}\right)}
\nonumber
\\
\nonumber
\label{equt:03_2a}
&=&    -\frac{\nu(1-s)}{(2-s)n}-
   \frac{(\nu^2/n+ \bar g(\cz_{i,k-1}(k/n))(1-s)}{2(2-s)^2 n}
-c_{\ref{equt:03_2}'}(1-s)n^{-3/2 }\\
\nonumber
&\geq& -\frac{(1-s)}{(2-s)n}\left(\nu+\frac{\Vert \bar g\Vert_{\infty}}{2}+\nu^2/n
 +2c_{\ref{equt:03_2}'}n^{-1/2 }\right)\\
\label{equt:03_4}
&\geq& -\frac{(1-s)}{(2-s)n}\left(\nu+\frac{\Vert \bar g\Vert_{\infty}}{2}+
 c_{\ref{equt:03_4}}n^{-1/2 }\right)
.\end{eqnarray}
Substituting in~(\ref{equt:03_2}) we get
\begin{eqnarray}
\nonumber
\EE\left( s^{\tilde{M}^n_k}\right)
\label{equt:03_5}
&\geq&  \EE\left( \left(\frac{1}{2-s}\right)^{\tilde{M}^n_{k-1}}
 \exp{\left\{-\tilde{M}^n_{k-1}
\frac{(1-s)}{(2-s)n}\left(\nu+\frac{\Vert \bar g\Vert_{\infty}}{2}+
 c_{\ref{equt:03_4}}n^{-1/2 }\right)\right\}}\right)\\
\nonumber
&\geq& \EE\left( \left(\frac{1}{2-s}
-\frac{(1-s)}{(2-s)^2n}\left(\nu+\frac{\Vert \bar g\Vert_{\infty}}{2}+
 c_{\ref{equt:03_4}}n^{-1/2 }\right)\right)^{\tilde{M}^n_{k-1}} \right)\\
&=:&  \EE\left(\tilde f(s)^{\tilde{M}^n_{k-1}} \right) .\end{eqnarray}
Let $f(s)$ be the generating function of the geometric distribution
with parameter $p=\frac12-\frac{b_n}{4n}$, then
\begin{eqnarray}
\nonumber
f(s)&=&\frac{1/2-b_n/4n}{1-(  1/2+b_n/4n)s}
\\  &=& \frac{1}{2-s}\left( 1-\frac{b_n(1-s)}{(2-s)n(1-\frac{b_ns}{2n(2-s)})}\right)
\nonumber
\\
&\leq& \frac{1}{2-s}\left( 1-\frac{b_n(1-s)}{(2-s)n}\right)
.\end{eqnarray}
If one takes $b_n=\nu+\frac{\Vert \bar g\Vert_{\infty}}{2}+
 c_{\ref{equt:03_4}}n^{-1/2 }$ then we get that
\[ \tilde f(s)\geq f(s),\;\; 0\leq s\leq 1,\]
and hence by iterating~(\ref{equt:03_5}) we get
\begin{eqnarray}
\nonumber
\EE\left( s^{\tilde{M}^n_k}\right)
&\geq&  \EE\left( s^{M^n_k}\right) ,\;\;0\leq s\leq 1.
\end{eqnarray}
Therefore
\begin{eqnarray}
\PP(\tilde{M}^n_{k} >0)\leq
\PP(M^n_{k}>0), \;\;\forall k\geq 1,
\end{eqnarray}
and hence by Lemma~\ref{lem:3}  we get that
\begin{eqnarray}
 \limsup_{n\rightarrow \infty} n \PP(\tilde{M}^n_{\lfloor n\delta\rfloor} >0)\leq
 \limsup_{n\rightarrow \infty} n \PP(M^n_{\lfloor n\delta\rfloor} >0)\leq  h(b,\delta)
.\end{eqnarray}
\gdm

\begin{lemma}
\label{lem:5} 
 Let $\tilde M^n, X^n$ be as above.
\begin{itemize}
 \item[{\bf (a)}]
For any $\delta>0$,
\[ \limsup_{n\rightarrow \infty} \EE\left( \tilde M^n_{\lfloor n\delta\rfloor}\right) \leq
 \limsup_{n\rightarrow \infty}\tilde M_0^n e^{b\delta}, \]
and hence,
\[ \limsup_{n\rightarrow \infty} \EE\left( X^n_{\lfloor n\delta\rfloor/n}(1)\right) \leq
\limsup_{n\rightarrow \infty} X^n_0(1) e^{b\delta}. \]
\item[{\bf (b)}]
For any $\delta, a>0$,
\[ \limsup_{n\rightarrow \infty} \PP\left( \sup_{k\leq \lfloor n\delta\rfloor}X^n_{k/n}(1)\geq a \right) \leq
 \frac{\limsup_{n\rightarrow \infty} X_0^n(1)\(e^{b\delta}\vee 1\)}{a}. \]
\end{itemize}
\end{lemma}
\paragraph{Proof:}{\bf (a)} The proof of goes along the similar lines as the proof
of the previous lemma. First  recall that
\begin{eqnarray}
\label{eq:28_12_1}
\EE\left( Z^n_{i,k}| X^n_{(k-1)/n}\right) &\leq& 1+ \frac{\nu}{n}+\frac{\Vert \bar g\Vert_{\infty}}{2n},
\end{eqnarray}
for all $i,k,n$. Hence, by iteration, we get
\begin{eqnarray}
 \EE\left(\tilde M^n_{\lfloor n\delta\rfloor}\right)&\leq&
\tilde M^n_0 \left(1+ \frac{\nu}{n}+\frac{\Vert \bar
g\Vert_{\infty}}{2n}\right)^{\lfloor n\delta\rfloor},
\end{eqnarray}
and the result follows.\\
{\bf (b)} For all $k\geq 0$, define 
\[ V^n_k= X^n_{k/n}(1) \left(1+ \frac{\nu}{n}+\frac{\Vert \bar
g\Vert_{\infty}}{2n}\right)^{-k}.\]
Then using~(\ref{eq:28_12_1}) it is easy to check that $\{V^n_k\}_{k\geq 0}$ is a nonnegative $\{\cF_k^{X^n}\}_{k\geq 0}$-supermartingale. Therefore by maximal inequalities for non-negative supermartingales we get
\begin{eqnarray}
\label{eq:28_12_2}
 \PP\left( \sup_{k\leq \lfloor n\delta\rfloor}V^n_k\geq a\right)&\leq&
 \frac{\EE\( V^n_{0}\)}{a}=\frac{X^n_0(1)}{a}.
\end{eqnarray}
To prove the result we consider the cases  $\nu + \frac{\Vert \bar
g\Vert_{\infty}}{2}<0$ and $\nu + \frac{\Vert \bar
g\Vert_{\infty}}{2}\geq 0$ separately. First suppose that  $\nu +
\frac{\Vert \bar g\Vert_{\infty}}{2}<0$. Recall the definition of
$V^n_k$ to get that, in this case,
\begin{eqnarray}
\label{eq:28_12_3}
 \lefteqn{\PP\left( \sup_{k\leq \lfloor n\delta\rfloor}V^n_k\geq a\right)}
 \\
 \nonumber&=&\PP\left( \sup_{k\leq \lfloor n\delta\rfloor}
 X^n_{k/n}(1) \left(1+ \frac{\nu}{n}+\frac{\Vert \bar
g\Vert_{\infty}}{2n}\right)^{-k}\geq a\right)\\
\nonumber
&\geq&\PP\left( \sup_{k\leq \lfloor n\delta\rfloor}X^n_{k/n}(1)\geq a\right).
\end{eqnarray}
By putting (\ref{eq:28_12_2}), (\ref{eq:28_12_3}) together we get that
\begin{eqnarray}
\label{eq:28_12_4}
\PP\left( \sup_{k\leq \lfloor n\delta\rfloor}X^n_{k/n}(1)\geq a\right)&\leq& \frac{X^n_0(1)}{a},\;\;{\rm for}\;
 \nu + \frac{\Vert \bar g\Vert_{\infty}}{2}<0.
\end{eqnarray}
Now let $\nu + \frac{\Vert \bar g\Vert_{\infty}}{2}\geq 0$. Then we get
\begin{eqnarray}
\label{eq:28_12_5}
 \lefteqn{\PP\left( \sup_{k\leq \lfloor n\delta\rfloor}V^n_k\geq a\left(1+ \frac{\nu}{n}+\frac{\Vert \bar
g\Vert_{\infty}}{2n}\right)^{-\lfloor n\delta\rfloor}\right)}
\\
\nonumber
&=&\PP\left( \sup_{k\leq \lfloor n\delta\rfloor}
 X^n_{k/n}(1) \left(1+ \frac{\nu}{n}+\frac{\Vert \bar
g\Vert_{\infty}}{2n}\right)^{-k}\geq a\left(1+ \frac{\nu}{n}+\frac{\Vert \bar
g\Vert_{\infty}}{2n}\right)^{-\lfloor n\delta\rfloor}\right)\\
\nonumber
&\geq&\PP\left( \sup_{k\leq \lfloor n\delta\rfloor}X^n_{k/n}(1)\geq a\right)
\end{eqnarray}
Apply this and~(\ref{eq:28_12_1}) with $a\left(1+ \frac{\nu}{n}+\frac{\Vert \bar
g\Vert_{\infty}}{2n}\right)^{-n}$ instead of $a$ to get
\begin{eqnarray}
\label{eq:28_12_6}
\lefteqn{\PP\left( \sup_{k\leq \lfloor n\delta\rfloor}X^n_{k/n}(1)\geq a\right)}\\
\nonumber
&\leq& \frac{X^n_0(1)\left(1+ \frac{\nu}{n}+\frac{\Vert \bar
g\Vert_{\infty}}{2n}\right)^{\lfloor n\delta\rfloor}}{a},\;\;{\rm for}\;
 \nu + \frac{\Vert \bar g\Vert_{\infty}}{2}\geq 0.
\end{eqnarray}
By combining (\ref{eq:28_12_4}), (\ref{eq:28_12_6}), we get
\begin{eqnarray}
\label{eq:28_12_7}
\PP\left( \sup_{k\leq \lfloor n\delta\rfloor}X^n_{k/n}(1)\geq a\right)&\leq& \frac{X^n_0(1)\(\left(1+ \frac{\nu}{n}+\frac{\Vert \bar
g\Vert_{\infty}}{2n}\right)^{\lfloor n\delta\rfloor}\vee 1\)}{a},
\end{eqnarray}
and by letting $n\rightarrow \infty$ we are done.
\gdm

The next result generalizes the previous lemma.
\begin{lemma}
\label{lem:nov1}
Let $f$ be a bounded non-negative measurable function  on $\RR^d$.
Then, for any $\delta>0$,
\[ \EE\left(  X^n_{\lfloor n\delta\rfloor/n}(f)\right) \leq
\left(1+ \frac{\nu}{n}+\frac{\Vert \bar
g\Vert_{\infty}}{2n}\right)^{\lfloor n\delta\rfloor}  X^n_0\left(S_{\delta} f\right), \]
where $\{S_t\}_{t\geq 0}$ is the semigroup of the Brownian motion.
\end{lemma}
\paragraph{Proof:} The proof goes along the similar lines as the proof
of the previous lemma. For any $k\geq 1$ we have 
\begin{eqnarray*}
\EE\left( X^n_{k/n}(f)| X^n_{(k-1)/n}\right)&=& \EE\left(\frac{1}{n}\sum_{i=1}^{\tilde
M^n_{k-1}}
  Z^n_{i,k} f(\cz_{i,k-1}(k/n)) | X^n_{(k-1)/n}\right)\\
&=&\frac{1}{n}\sum_{i=1}^{\tilde M^n_{k-1}} \EE\left(\left. \EE\left(\left.
  Z^n_{i,k} f(\cz_{i,k-1}(k/n))\right|
  X^n_{\frac{k}{n}-} \right) \right| X^n_{(k-1)/n}\right)
\\
&=& \frac{1}{n}\sum_{i=1}^{\tilde M^n_{k-1}}\EE\left(\left. f(\cz_{i,k-1}(k/n))
\EE\left(\left.
  Z^n_{i,k} \right| X^n_{\frac{k}{n}-} \right)\right| X^n_{(k-1)/n}\right)
\\
&\leq&\left(1+ \frac{\nu}{n}+\frac{\Vert \bar
g\Vert_{\infty}}{2n}\right) \frac{1}{n} \sum_{i=1}^{\tilde M^n_{k-1}}
\EE\left(\left. f(\cz_{i,k-1}(k/n)) \right| X^n_{(k-1)/n}\right)
\\
&=&\left(1+ \frac{\nu}{n}+\frac{\Vert \bar
g\Vert_{\infty}}{2n}\right) \frac{1}{n} \sum_{i=1}^{\tilde M^n_{k-1}}
S_{1/n}f(\cz_{i,k-1}((k-1)/n))
\\
&=&\left(1+ \frac{\nu}{n}+\frac{\Vert \bar
g\Vert_{\infty}}{2n}\right) X^n_{(k-1)/n}\left( S_{1/n}f\right),
\end{eqnarray*}
for all $k,n$ and hence by iteration
\begin{eqnarray}
\EE\left(  X^n_{\lfloor n\delta\rfloor/n}(f)\right)&\leq&
\left(1+ \frac{\nu}{n}+\frac{\Vert \bar
g\Vert_{\infty}}{2n}\right)^{\lfloor n\delta\rfloor}
X^n_0\left(S_{\lfloor n\delta\rfloor/n} f\right),
\end{eqnarray}
and the result follows.
\gdm

\section{Tightness of the contour process}
\label{sec-tcont}
\setcounter{equation}{0}
In this section we will prove the tightness of the sequence of the contour processes $\{Y^n\}_{n\geq 1}$.
The 
following proposition is the main result of
this section.
\begin{prop}[Tightness of $\{Y^n\}_{n\geq 1}$]
\label{prop:tight}
For any $\delta>0$, $T>0$,
\begin{eqnarray}
\lim_{\ep\downarrow 0} \limsup_{n\rightarrow\infty}\PP\left(\sup_{0\leq 
t\leq T}\sup_{0\leq s\leq  \ep}|Y^n_{t+s}-Y^n_{t}|>\delta\right)=0,
\end{eqnarray}
that is, $\{Y^n\}_{n\geq 1}$ is $C$-tight in $D_{\RR}[0,\infty)$.
\end{prop}
The proof of the proposition will be given in this section. Recall
the definition of the discrete  version of the local time for $Y^n$
and its inverse (see 
(\ref{eq:1.3}), (\ref{eq:30_12_10}) and
(\ref{eq:1.4}) for the same definitions in the case without
environment). Fix an arbitrary $c_0>0$. We will first handle the
tightness on the time interval
\[ t\in [0, \tau^{n,0}_{c_0}],\]
and we start with the following proposition.
\begin{prop}[Tightness of $\{Y^n\}_{n\geq 1}$ --- no jumps up]
\label{prop:1}
For any $\delta>0$,
\begin{eqnarray}
\lim_{\ep\downarrow 0} \limsup_{n\rightarrow\infty}\PP\left(\sup_{0\leq t\leq  \tau^{n,0}_{c_0}}\sup_{0\leq s\leq  \ep}(Y^n_{t+s}-Y^n_{t})_+>\delta\right)=0.
\end{eqnarray}
\end{prop}
The proof of the proposition will be given
after we present several preliminary lemmas. For any $a\geq 0$ and
$r_1<r_2$, recall the measure-valued  process
$X^{n,r_1,r_2}_{a_n,a_n+t}$, see (\ref{eq:31_12_01}). Fix an
arbitrary $\delta>0$. 
Recall that $b$ was defined in~(\ref{equt:nov6}).
We have the following lemma.
\begin{lemma}\label{lem3b2} For any $\delta>0, r>0, \ep'>0,$
\begin{eqnarray}
\limsup_{n\rightarrow \infty}
\PP\left(X^{n,r,r+\ep'}_{a,a_{n}+\delta}(1)>0\right)&\leq &
\ep'h(b,\delta).
\end{eqnarray}
\end{lemma}
\paragraph{Proof:}
\begin{eqnarray}
 \PP\left(X^{n,r,r+\ep'}_{a,a_{n}+\delta}(1)>0\right)&\leq &
\sum_{i=0}^{\lfloor n\ep'\rfloor}\PP\left(X^{n,r+i/n,r+(i+1)/n}_{a,a_{n}+\delta}(1)>0\right).
\end{eqnarray}
Since, by Lemma~\ref{lem:4},
\[ \limsup_{n\rightarrow \infty}n \PP\left(X^{n,r+i/n,r+(i+1)/n}_{a,a_{n}+\delta}(1)>0\right)
\leq h(b,\delta),\]
 the result follows.
\gdm

The following corollary is immediate.
\begin{cor}
\label{cor:2} For any $\delta>0, r>0, \ep'>0,$ 
\begin{eqnarray}
\label{cor:1}
\limsup_{n\rightarrow \infty} \PP(\sup_{\tau^{n,a}_r\leq t\leq \tau^{n,a}_{r+\ep'}} (Y^n_t-Y^n_{\tau^{n,a}_r})_+ >\delta)&\leq&\ep' h(b,\delta).
\end{eqnarray}
\end{cor}

The next corollary gives a bound on the positive increment of
$Y^n$.
\begin{cor}
\label{cor:3}
For any $\delta>0, r>0,$
\begin{eqnarray}
\lim_{\ep\downarrow 0}\limsup_{n\rightarrow \infty} \PP(\sup_{\tau^{n,a}_r\leq t\leq \tau^{n,a}_r+\ep} (Y^n_t-Y^n_{\tau^{n,a}_r})_+ >\delta)&=&0.
\end{eqnarray}
\end{cor}
\paragraph{Proof:}
We first
prove that  for any  $\ep'>0$,
\begin{eqnarray}
\label{equt:2}
\lim_{\ep\downarrow 0}\limsup_{n\rightarrow\infty}\PP\left(\tau^{n,a}_{r+\ep'}\leq \tau^{n,a}_r+\ep\right)=0.
\end{eqnarray}
Suppose that there exist deterministic $\ep'>0,\delta''\in (0,1/2)$
and subsequences
$n_k\rightarrow \infty, \ep_k\downarrow 0$ such that
\begin{eqnarray}
\lim_{k\rightarrow\infty}
\PP\left(\tau^{n_k,a}_{r+\ep'}\leq \tau^{n_k,a}_r+\ep_k\right)\geq \delta''.
\end{eqnarray}
To avoid cumbersome notation,
for the rest of the proof we write $n$ and $\ep_n$ for
$n_k$ and $\ep_k$ respectively.
Note that it follows from~(\ref{equt:2}) that
\begin{eqnarray}
 \ell^{n,a_{n}}_{\tau^{n,a}_{r+\ep'}}-\ell^{n,a_{n}}_{\tau^{n,a}_r}=\ep'.
\end{eqnarray}
Then as in Lemma \ref{lem3b2}, we may define the sequence of measure-valued processes $X^{n}$ with total mass
\begin{eqnarray}
X^{n}_{a,a_{n}+s}(1)&=&
\ell^{n,a_{n}+s}_{\tau^{n,a}_{r+\ep'}}-\ell^{n,a_{n}+s}_{\tau^{n,a}_r}.
\end{eqnarray}
This process starts at the total mass
\[X^{n}_{a,a_{n}}(1)=\ep',\]
and, appealing to \cite{M}, as $n\rightarrow \infty$, it converges
weakly in $D_{\RR}[0,\infty)$ to the continuous process $s\mapsto
X_{a,a+s}(1)$ starting at $\ep'$. Therefore, by the weak convergence
properties, there exists $\bar\delta>0$ such that 
%
\begin{eqnarray}\label{eq3b13}
\liminf_{n\rightarrow\infty}\PP\left( \inf_{s\leq \bar\delta}\(
\ell^{n,a_{n}+s}_{\tau^{n,a}_{r+\ep'}}-\ell^{n,a_{n}+s}_{\tau^{n,a}_r}\)\geq
\ep'/2\right)\geq 1-\delta''.
\end{eqnarray}
Note that on the event in (\ref{eq3b13}), we have 
\[\frac{\ep'}{2}\le \ell^{n,a_{n}+s}_{\tau^{n,a}_{r+\ep'}}-\ell^{n,a_{n}+s}_{\tau^{n,a}_r}
\le n^{-1}\sum_{l:\;\tau^{n,a}_r< l/n^2\leq
\tau^{n,a}_{r+\ep'}} 1_{Y^{n}_{l/n^2}= a_{n}+s}.\] Summing
over $s=\frac{i}{n}\le\bar{\de}$, we get that with probability
greater than $1-2\de''>0 $, the occupation time
 of $Y^{n}$ on the time interval
 $(\tau^{n,a}_r, \tau^{n,a}_{r+\ep'}]$ is bounded from below by
\begin{equation}
  \frac{1}{n^2}\sum_{l:\;\tau^{n,a}_r< l/n^2\leq
\tau^{n,a}_{r+\ep'}} 1_{Y^{n}_{l/n^2}\geq a_{n}} \label{equt:1}
\geq\frac{1}{2}\bar\delta \ep'>0.
\end{equation}
On the other hand the total occupation time of $Y^{n}$ on the interval
\[[\tau^{n,a}_r,\tau^{n,a}_{r+\ep'}]\subset[\tau^{n,a}_r,\tau^{n,a}_r+\ep_n] \]
is bounded by
$ \ep_n\downarrow 0$,
which contradicts~(\ref{equt:1}). Hence~(\ref{equt:2}) follows.

Continuing with the proof of the lemma, we have from (\ref{equt:2}) that
for any $\ep'>0$,
\begin{eqnarray}
\lefteqn{\lim_{\ep\downarrow 0}\limsup_{n\rightarrow\infty}
\PP\left(\sup_{\tau^{n,a}_r\leq t\leq \tau^{n,a}_r+\ep}(Y^n_t-Y^n_{\tau^{n,a}_r})_+\geq \delta\right)}
\\
\nonumber
&\leq& \lim_{\ep\downarrow 0}\limsup_{n\rightarrow\infty}
\PP\left(\sup_{\tau^{n,a}_r\leq t\leq \tau^{n,a}_{r+\ep'}}(Y^n_t-Y^n_{\tau^{n,a}_r})_+\geq \delta\right)\\
\nonumber
&&\mbox{}+
\lim_{\ep\downarrow 0}\limsup_{n\rightarrow\infty}\PP\left(\tau^{n,a}_{r+\ep'}\leq \tau^{n,a}_r+\ep\right)
\leq \ep' h(b,\delta),
\end{eqnarray}
where the last inequality follows by~(\ref{equt:2}) and Corollary~\ref{cor:2}.
Since $\ep'$ was arbitrary we are done.
\qed

We now introduce further notation. Let 
\begin{eqnarray}
\label{01_06}
\bar{\ell}^{n,a}&\equiv& \ell^{n,a}_{\tau^{n,0}_{c_0}}.
\end{eqnarray}
We 
will prove the
following lemma.
\begin{lemma}
    For any $\delta>0$,
\label{lem:2}
\begin{eqnarray}
\lim_{\ep\downarrow 0} \limsup_{n\rightarrow\infty}\PP\left(\sup_{0\leq r\leq \bar{\ell}^{n,a}}
\sup_{0\leq s\leq\ep}
(Y^n_{\tau^{n,a}_r+s}-Y^n_{\tau^{n,a}_r})_+>\delta\right)=0.
\end{eqnarray}
\end{lemma}
\paragraph{Proof:}
 Now we will need some further notation.
Denote
\begin{eqnarray*}
N^{a,n}_{s}&=&\{\mbox{number of excursions of $Y^n$ starting at  level $a_n$ }\\
&&\mbox{above the  level $a_n+s$ on the time interval $[0,\tau^{n,0}_{c_0}]$}\}\\
&=&\{\mbox{number of particles in the original branching particle system  }\\
&&\mbox{at time $a_nn$ whose descendants survive till time
$(a_n+s)n$}\}
\end{eqnarray*}
By Lemmas~\ref{lem:4}, \ref{lem:5} and the Markov property of the branching system
we immediately get
\begin{eqnarray}
\nonumber
\limsup_{n\rightarrow \infty} \EE\left( N^{a,n}_{s}\right) &\leq& \limsup_{n\rightarrow \infty}
 \EE\left( X^n_a(1)\right) n  \PP(\tilde{M}^n_{\lfloor n s\rfloor} >0|\tilde{M}^n_0=1)
\\
\label{equt:03_6}
&\leq&c_0 e^{ba}h(b,s).
\end{eqnarray}
For $i=1,\ldots N^{a,n}_{\delta/2}$ define 
\[ \sigma^n_i=\inf\{t>\hat{\tau}^n_{i-1}:\; Y^n_t \geq a_n+\delta/2\},\]
where
\begin{eqnarray*}
\hat{\tau}^n_0&=& 0,\\
\hat{\tau}^n_{i}&=& \inf\{t>\sigma^n_{i}:\; Y^n_t = a_n\}.
\end{eqnarray*}
That is, $ \sigma^n_i$ are the times when successful excursions of
$Y^n$ reach the level $a_n+\delta/2$. Then we have, for any fixed integer
$m$,
\begin{eqnarray}
\label{equt:3}
\lefteqn{ \lim_{\ep\downarrow 0}\limsup_{n\rightarrow\infty}\PP\left(\sup_{0\leq r\leq  \bar{\ell}^{n,a}}
\sup_{0\leq s\leq\ep}
(Y^n_{\tau^{n,a}_r+s}-Y^n_{\tau^{n,a}_r})_+>\delta\right)}\\
\nonumber
&\leq& \lim_{\ep\downarrow 0} \limsup_{n\rightarrow\infty}\EE\left(\sum_{i=1}^{N^{a,n}_{\delta/2}}
 \PP\left(\sup_{0\leq s\leq\ep}(Y^n_{\sigma^n_i+s}-Y^n_{\sigma^n_i})_+>\delta/2|
 X^n_{a_n+\delta/2}\right); N^{a,n}_{\delta/2}\leq m
\right)\\
\nonumber
&&\mbox{}+ \lim_{\ep\downarrow 0}\limsup_{n\rightarrow\infty}\PP\left(  N^{a,n}_{\delta/2}> m\right).
\end{eqnarray}
By an argument similar to the one in Corollary~\ref{cor:3} we get that
\begin{eqnarray}
\label{equt:6}
\lim_{\ep\downarrow 0}\limsup_{n\rightarrow \infty}
\PP\left(\sup_{0\leq s\leq\ep}(Y^n_{\sigma^n_i+s}-Y^n_{\sigma^n_i})_+>\delta/2|
 X^n_{a_n+\delta/2}\right)
=0.
\end{eqnarray}
This implies
that
\begin{eqnarray*}
&&
\!\!\!\!\!\!\!
\!\!\!\!\!\!\!
\!\!\!\!\!\!\!
\!\!\!\!\!\!\!
\!\!\!\!\!\!\!
\!\!\!\!\!\!\!
\!\!\!\!\!\!\!
\lim_{\ep\downarrow 0}\limsup_{n\rightarrow\infty}
\PP\left(\sup_{0\leq r\leq  \bar{\ell}^{n,a}}
\sup_{0\leq s\leq\ep}
(Y^n_{\tau^{n,a}_r+s}-Y^n_{\tau^{n,a}_r})_+>\delta\right)
\\
&\leq&  \lim_{\ep\downarrow 0}\limsup_{n\rightarrow\infty}\PP\left(  N^{a,n}_{\delta/2}>  m\right)
\leq \frac{c_0 e^{ba}h(b,\delta/2)}{m},
\end{eqnarray*}
where the last inequality follows by the Markov inequality and
(\ref{equt:03_6}). Since $m$ was arbitrary we are done.
\qed

We can now complete the proof of
Proposition~\ref{prop:1}.
\paragraph{Proof of Proposition~\ref{prop:1}:}
For $\delta>0$
let
\begin{eqnarray}
\cT^{i,\delta}_n&=&\{t\leq  \tau^{n,0}_{c_0}:\;
Y^n_t \in [i\delta/2, (i+1)\delta/2]\}.
\end{eqnarray}
Then
\begin{eqnarray}
\nonumber
\lefteqn{\PP\left(\sup_{0\leq t\leq  \tau^{n,0}_{c_0}}\sup_{0\leq s\leq  \ep}
(Y^n_{t+s}-Y^n_{t})_+>\delta\right)}\\
\nonumber
&\leq&
\sum_{i=0}^{\lfloor 2K_1/\delta\rfloor}
\PP\left(\sup_{t\in \cT^{i,\delta}_n}\sup_{0\leq s\leq  \ep}
(Y^n_{t+s}-Y^n_{t})_+>\delta
\right)
\\
\label{equt:8}
&\leq&
\sum_{i=1}^{\lfloor 2K_1/\delta\rfloor +1}
\PP\left(\sup_{r\leq  \bar{\ell}^{n,i\delta/2}}\sup_{0\leq s\leq  \ep}
(Y^n_{\tau^{n,i\delta/2}_r+s}-Y^n_{\tau^{n,i\delta/2}_r})_+>\delta/2
\right).
\end{eqnarray}
However by Lemma~\ref{lem:2} we get that, for every $i$,
\[ \lim_{\ep\rightarrow 0}\limsup_{n\rightarrow \infty}
\PP\left(\sup_{r\leq  \bar{\ell}^{n,i\delta/2}}\sup_{0\leq s\leq  \ep}
(Y^n_{\tau^{n,i\delta/2}_r+s}-Y^n_{\tau^{n,i\delta/2}_r})_+>\delta/2
\right)=0,\]
and this finishes the proof of Proposition~\ref{prop:1}.
\qed

To handle downward jumps, we need the following proposition.
\begin{prop}[Tightness of $\{Y^n\}_{n\geq 1}$ --- no jumps down]
\label{prop:2}
For any $\delta>0$,
\begin{equation}
\label{of-220209}
\lim_{\ep\downarrow 0} \limsup_{n\rightarrow\infty}\PP\left(\sup_{0\leq t\leq  \tau^{n,0}_{c_0}}\sup_{0\leq s\leq  \ep}(Y^n_{t+s}-Y^n_{t})_->\delta\right)=0.
\end{equation}
\end{prop}
\paragraph{Proof:} In fact the proof is easy if one considers
the process $Y^n$ 
reversed in time, that is the process $\tilde
Y_t^n=Y^n_{\tau^{n,0}_{c_0}-t}$, which is easily seen (see the
explicit argument below) to possess the same law (with $0\leq t\leq
\tau^{n,0}_{c_0}$) as the original process $Y^n_\cdot$. 
Since any jump down for $Y^n$ becomes a jump up
for $\tilde Y^n$, the claim (\ref{of-220209}) follows from
Proposition \ref{prop:1} applied to $\tilde Y^n$.

To see the reversibility claim, we introduce a sequence of
path transformations $\{T_z\}_{z=0,1/n,\ldots,K_1-1/n}$
on $\{(\bW^n_{t},Y^n_t)\}_{t=0,1/n^2,\ldots,\tau^{n,0}_{c_0}}$, each of
which is  measure preserving and preserves $\tau^{n,0}_{c_0}$, such
that
$$\{(\tilde \bW^n_{t},\tilde Y_t^n)\}_{t=0,1/n^2,\ldots,\tau^{n,0}_{c_0}} =
T_{K_1-1/n}\circ T_{K_1-2/n}\circ\ldots\circ T_0
\{(\bW^n_{t},Y_t^n)\}_{t=0,1/n^2,\ldots,\tau^{n,0}_{c_0}}\,,$$ where
$(\tilde \bW^n_{\cdot},\tilde Y_\cdot^n)$ denotes the image of
$(\bW^n_{\cdot},Y^n_\cdot)$
under the transformations. This will prove the claim.

To avoid cumbersome notation, we consider the case of $n=1$ only,
and we omit the index $n$. The general $n$ can be treated the same
way with proper scaling. For $z=0$, the transformation $T_0$ is
obtained as follows. If $t=\tau^0_j$ for some integer $j\in
\{0,1,\ldots,c_0\}$, that is $t$ is a return time of $Y_\cdot$ to
$0$, then define $t'(t)= \tau^0_{c_0}-\tau^0_j$. If $t\in
(\tau^0_{j-1},\tau^0_j)$ for some $j\in \{1,\ldots,c_0\}$, that is
$t$ belongs to the $j$th excursion of $Y_\cdot$ from $0$, then
define $t'(t)=t'(\tau^0_j)+t-\tau^0_{j-1}$. Then,
$$T_0(\bW_{\cdot},Y_\cdot)(t)
= (\bW_{t'(t)},Y_{t'(t)}).$$ (In words, $T_0$ reverses the order of
the excursions from $0$ but keeps the time orientation of each
excursion intact; Thus, the total length of the excursions is
preserved.) It is straightforward to check that the law of
$T_0(\bW_{\cdot},Y_\cdot)$ is the same as that of
$(\bW_{\cdot},Y_\cdot)$.

For $z=1$, let
$$t_z(1)=\min\{k> 0: Y_k=z\}, s_z(1)=\min\{k>t_z(1): Y_{k}=z,
Y_{k+1}=z-1\},$$ and for $j\geq 1$,
\begin{eqnarray*}
t_z(j+1)&=&\min\{k>  s_z(j): Y_{k}=z\},\\
 s_z(j+1)&=&\min\{k> t_z(j+1): Y_{k}=z,
Y_{k+1}=z-1\}.
\end{eqnarray*}
 $T_1$ is then defined as $T_0$ applied to the excursions
of the path {\it from level $z$}. Explicitly, let $j_z=\max\{j:
t_z(j)<\tau^0_{c_0}\}$. For $t\in \left(\cup_{j=1}^{j_z}
[t_z(j),s_z(j)]\right)^c$, set $t'(t)=t$. For each $j$, let
$\bar{t}_z(j,0)=t_z(j)$,
$\bar{t}_z(j,\ell)=\min\{t>\bar{t}_z(j,\ell):\;Y_t=z\}$, and
$\ell_z(j)=\max\{\ell:\;\bar{t}_z(j,\ell)=s_z(j)\}$. Let $t'$ be
defined on the interval $[t_z(j),s_z(j))$ in the same way as the case
of $z=0$ with $\tau^0_j,\;j=0,1,\cdots,c_0$ replaced by
$\bar{t}_z(j,\ell),\;\ell=0,1,\cdots,\ell_z(j)$. Then,
$$T_z(\bW_{\cdot},Y_\cdot)(t)
= (\bW_{t'(t)},Y_{t'(t)}).$$ Again, in words,
 $T_z$ reverses the order of the excursions from $z$ but keeps
the time orientation of each excursion intact; Thus, the total
length of the excursions is preserved.) It is straightforward to
check that the law of $T_z(\bW_{\cdot},Y_\cdot)$ is the same as that
of $(\bW_{\cdot},Y_\cdot)$. We can continue this procedure for
$z=2,3,\cdots,K_1-1$. As explained above, this completes the proof.
\qed

To 
finish the
proof of the Proposition~\ref{prop:tight} we need the folowing
lemmas that describe the limiting behavior of
$\{\bar\ell^{n,\cdot}\}_{n\geq 1}$ and $\{\tau^{n,0}\}_{n\geq 1}$
(recall that $\bar\ell^{n,\cdot}=
 \ell^{n,\cdot}_{\tau^{n,0}_{c_0}}$
was introduced in~(\ref{01_06})). 
\begin{lemma}
\label{lem:01_06_1}
For any $c_0>0$, the sequence of processes $\{\bar \ell^{n,\cdot}\}_{n\geq 1}$ is $C$-tight in $D_{\RR}$.
\end{lemma}
\paragraph{Proof:}
First recall from~(\ref{01_06}) and (\ref{eq:Oc5_1}),   that $X^{n,c_0}_{0,a}(1)=\bar \ell^{n,a}$ is
the total mass at time $a$ of the measure-valued  process $X^{n,c_0}_{0,\cdot}$
defined in the introduction. Since the sequence of measure-valued
processes $\{X^{n,c_0}_{0,\cdot}\}_{n\geq 1}$ is $C$-tight in $D_{\cM_F}$ (see \cite{M} and the
comments leading to (\ref{MP-3})), we get the desired result.
\qed

The next lemma studies the limiting
behavior of $\{\tau^{n,0}\}_{n\geq 1}$. Toward this end, recall that
according to our conventions introduced after~(\ref{eq:30_12_10}),
we use the same notation for an increasing function and the
corresponding measure. 
\begin{lemma}
\label{lem:01_06_2}
\begin{itemize}
 \item[{\bf (a)}] For any $r>0$,
 the sequence of random variables $
\{\tau^{n,0}_r\}_{n\geq 1}$ is tight and any limit point $\tau^0_r$ satisfies
\[ \PP(\tau^0_r=0)=0.\]
\item[{\bf (b)}] For any $\ep>0$,
    $A>0$, there exists $R>0$, such that
\[ \liminf_{n\rightarrow \infty}\PP\left((\tau^{n,0}_R>A\right)\geq 1-\ep. \]
 \item[{\bf (c)}]
 The sequence $
\{\tau^{n,0}\}_{n\geq 1}$ is tight in $\cM(\RR_+)$.
\item[{\bf (d)}] Let $\tau^0\in \cM(\RR_+)$ be an arbitrary limit point of
  $\{\tau^{n,0}\}_{n\geq 1}$. Then for any fixed $r\in \RR_+$, $\tau^0_t$ 
is continuous at $t=r$ with probability $1$.
\end{itemize}
\end{lemma}
\paragraph{Proof:}{\bf (a)}
Define
\[T^n_t(y)=\int^y_0 \ell^{n,z}_t dz=n^{-2}\sum^{\lfloor n^2 t\rfloor}_{i=0}
{\mathbf 1}_{Y^n_{n^{-2}i}\le
y}.\]
Note that
\begin{equation}
\label{01_06_2}
\tau^{n,0}_r = T^{n}_{\tau^{n,0}_r}(K_1).
\end{equation}
On the other hand
$$ T^{n}_{\tau^{n,0}_r}(K_1)= \int^{K_1}_0 \ell^{n,z}_{\tau^{n,0}_r} dz
\leq K_1 \sup_{s\leq K_1} \ell^{n,s}_{\tau^{n,0}_r},
$$
and since by Lemma~\ref{lem:01_06_1}, $\{\ell^{n,\cdot}_{\tau^{n,0}_r}\}_{n\geq 1}$ is tight, by~(\ref{01_06_2}) we get the tightness of $\{\tau^{n,0}_r\}_{n\geq 1}$.

Similarly, since $\{\ell^{n,\cdot}_{\tau^{n,0}_r}\}_{n\geq 1}$ is $C$-tight for any $\ep>0$ we can fix $\delta$ such that
\[  \PP(\inf_{s\leq \delta} \ell^{n,s}_{\tau^{n,0}_r}\geq c_0/2)\geq 1-\epsilon\]
for all $n$ sufficiently large.
Using this, ~(\ref{01_06_2}) and 
the definition of $T^n$ we get
\[
\tau^{n,0}_r = T^{n}_{\tau^{n,0}_r}(K_1) \ge \int_0^{\delta}
\ell^{n,s}_{\tau^{n,0}_r} ds \geq \frac{r}{2}\delta
\]
with probability at least $1-\epsilon$ for all $n$ sufficiently
large. Since $\epsilon$ was arbitrary we get that any limit point of
$\tau^{n,0}_r$ is greater than $0$ with probability $1$.
\\
{\bf (b)} For any $K>0$ we can represent
\[ \tau^{n,0}_{Kr}=\sum_{i=1}^{K} \tau^{n,0}_{i,r},\]
where, for each $i$, $\tau^{n,0}_{i,r}$ is distributed as $\tau^{n,0}_r$.
Fix arbitrary $\ep, A>0$. Since, by part {\bf (a)} of the lemma,
any limit point of $\tau^{n,0}_{i,r}$ is strictly greater than
$0$ with probability one, we can easily choose $K$
sufficiently large such that
$ \tau^{n,0}_{K r}=\sum_{i=1}^{K} \tau^{n,0}_{i,r}>A$ with probability at least $1-\ep$, for all $n$ sufficiently large.
\\
{\bf (c)} Immediate from {\bf (a)}.
\\
{\bf (d)} Let $\tau^0\in \cM(\RR_+)$ be a limiting point
 $\{\tau^{n,0}\}_{n\geq 1}$.
 To prove this part of the lemma we have to show that, for any $\ep>0$, there exists $\delta>0$, such that
\begin{eqnarray}
\label{eq:28_12_10}
\PP\(\tau^0_{r+\delta}-\tau^0_{r-\delta}>\ep\)\leq \ep.
\end{eqnarray}
Similarly to what we have done in {\bf (a)} define, 
\begin{eqnarray}
T^{n,r}_{s,t}(y)&=&\int^y_0 (\ell^{n,z}_t - \ell^{n,z}_s)
dz,\;\; 0\leq s\leq t.
\end{eqnarray}
Then we have
\begin{eqnarray}
\label{eq:28_12_9}
\tau^{n,0}_{r+\delta}-\tau^{n,0}_{r-\delta}&=&
T^{n,r}_{\tau^{n,0}_{r-\delta},\tau^{n,0}_{r+\delta}}(K_1)
=\int^{K_1}_0
(\ell^{n,z}_{\tau^{n,0}_{r+\delta}} - \ell^{n,z}_{\tau^{n,0}_{r-\delta}}) dz\\
&=& \int^{K_1}_0 X^{n,r-\delta, r+\delta}_{s}(1) ds \leq K_1
\sup_{s\leq K_1} X^{n,r-\delta, r+\delta}_{0,s}(1),
\nonumber\end{eqnarray} 
where
recall that $X^{n,r-\delta, r+\delta}_{0,s}$ is the measure-valued
process corresponding to the branching particle system in random
environment, constructed in Section~\ref{sec:1}
(see~(\ref{eq:31_12_01})), that starts at time $s=0$ with initial
mass $2\delta$. By Lemma~\ref{lem:5}(b)
\begin{eqnarray}
\label{eq:28_12_8}
\PP\(\sup_{s\leq K_1} X_s^{n,r-\delta, r+\delta}(1)> \ep\)&\leq& \frac{2 X_0^{n,r-\delta, r+\delta}(1) (e^{bK_1}\vee 1)}{\ep}\\
\nonumber
 &=&\frac{4\delta (e^{bK_1}\vee 1)}{\ep},
\end{eqnarray}
for all $n$ sufficiently large. We can take $\delta$ sufficiently small such that the right hand side of~(\ref{eq:28_12_8}) is less than $\ep/2$, and this together with~(\ref{eq:28_12_9}) implies that
\begin{eqnarray}
\PP\(\tau^{n,0}_{r+\delta}-\tau^{n,0}_{r-\delta}>\ep\)\leq \ep/2.
\end{eqnarray}
for all $n$ sufficiently large. Therefore~(\ref{eq:28_12_10}) follows for any limit point of $\{\tau^{n,0}\}_{n\geq 1}$.
\qed

Now we are ready to complete the proof of Proposition~\ref{prop:tight}.
\paragraph{Proof of Proposition~\ref{prop:tight}:} Proposition~\ref{prop:tight} follows immediately from Propositions~\ref{prop:1}, \ref{prop:2}, Lemma~\ref{lem:01_06_2}{\bf (b)}, and the fact that $c_0$ was arbitrary.
\qed

\section{Tightness of $\{(\WW^n,\ell^n)\}_{n\geq 1}$ and proof of Theorem~\ref{thm:1}
}
\label{sec-4}

\setcounter{equation}{0}

The bulk of this section is devoted to the proof of the following
proposition. 
\begin{prop}
\label{prop:3}
The sequence $\{(\WW^n,\ell^n, \tau^{n,0})\}_{n\geq 1}$ is tight in
\mbox{$D_{\cW\times\cM(\RR_+)}$} $\times$
\mbox{$\cM(\RR_+)$.} Let  $(\WW,\ell,\tau^0)$ be its arbitrary limiting point.
Then $(\WW,\ell,\tau^0)$  belongs to
$C_{\cW\times\cM(\RR_+)}\times \cM(\RR_+)$.
Moreover,
$\ell$ is the local time of $Y$ ($Y$ is the lifetime
of $\WW$),
that is,
\begin{equation}
    \label{limpointloctime}
    \int_0^t {\mbox{\bf 1}}_{Y_s\leq a}\, ds =
\int_0^a \ell^r_t\,dr,\;\;\forall a\geq 0, \; t\geq 0.
\end{equation}
\end{prop}
Note that following our conventions,
we denote by $\ell^r(dt)$ the measure and by
$\ell^r_t=l^r([0,t])$ the corresponding increasing distribution function
corresponding to $\ell$.

The proof of Proposition \ref{prop:3} is long and we indicate the
main steps. We will first prove the tightness of the sequence of
processes $\{\WW^n\}_{n\geq 1}$, based on the tightness of the
contour process established in Section \ref{sec-tcont}. This will be
obtained in Lemma \ref{lem:07_06_1}, after going through a fair
amount of preliminary material. The tightness of the sequence of the
local time process $\{\ell^n\}_{n\geq 1}$ is then obtained in Lemma
\ref{lem-newtl}, thus completing the proof of Proposition
\ref{prop:3}. The rest of the section is devoted to the
identification of the limiting snake representation. 
Here we have to  identify a limit point of the sequence of the local times $\{\ell^n\}_{n\geq 1}$
as the local time of the limiting contour process, and this is done in Lemma
\ref{lem-loctime}. Additionally,  in Lemma~\ref{lem:28_12_1} we verify that
 a limiting point of $\{\tau^{n,0}_{c_0}\}_{n\geq 1}$
is  indeed the value at $c_0$ of the inverse function of the limiting local time.
 The proof of Theorem~\ref{thm:1} is an immediate corollary
of these facts, and is presented at the end of the section.

As in the previous section, where the tightness of the contour
processes $\{Y^n\}_{n\geq 1}$ was obtained, we 
first
handle tightness on the time interval $[0,\tau^{n,0}_{c_0}]$. Fix an
arbitrary $a\in [0,K_1)$ and recall that $\{X^{n,c_0}_{0,t}\}_{t\geq
0}$ (see~(\ref{eq:Oc5_1})) is the measure-valued process
characterising the branching particle picture, and in particular,
$nX^{n,c_0}_{0,a_n}(1)$ is the number of particles alive at time
$a_n={\lfloor an\rfloor}/{n}$. First we derive a bound on the
maximal displacement of the offsprings from the ancestors during the
time interval $[a_n\,, a_n+\delta]$. This estimate will be crucial
for proving tightness of paths of the Brownian snake in random
environment.

Fix $\eta\in(0,1/4)$ arbitrary small. Define
\begin{eqnarray}
    Z^{n,\eta}_{a_n\,,\delta}&=& n\int_0^{\tau^{n,0}_{c_0}}\mbox{\bf 1}_
    {   (|\hat{\WW}^n_s-
    \bW^n_s(a_n)|>\delta^{1/2-\eta})}\ell^{a_n+\delta}(ds)\\
\nonumber
&=& \#   \{\mbox{particles alive
at time $a_n+\delta$ that are displaced by more than}
\\
\nonumber
&&\mbox{ \ \ \
$\delta^{1/2-\eta}$
from the
 ancestor at time $a_n$}\}.
\end{eqnarray}
\begin{lemma}
\label{lem:nov4} There exists
$\delta_{\ref{lem:nov4}}>0$, such that
\begin{eqnarray}
\PP(Z^{n,3\eta/2}_{a_n\,,\delta}>0)\leq e^{-\delta^{-\eta}},\;\; \forall \delta\leq \delta_{\ref{lem:nov4}}\,, \forall n.
\end{eqnarray}
\end{lemma}

We postpone the proof of Lemma \ref{lem:nov4}, and prepare
some preliminary estimates.
Introduce the event
$$ {\cal W}_{n,a,\delta,k,\eta,s}=
\{
|\bW^n_s(a_n+\delta-\delta2^{-k})
-\bW^n_s(a_n+\delta-\delta 2^{-(k-1)})|>\delta^{1/2-\eta}2^{-k/4}\}\,,$$
and define
$$\tZ^k_{a_n\,,\delta} =
n\int_0^{\tau^{n,0}_{c_0}}
\mbox{\bf 1}_{ {\cal W}_{n,a,\delta,k,\eta,s}}
\ell^{n,a_n+\delta-\delta 2^{-k-1}}(ds)\,,$$
which gives the number of particles alive at time $a_n+\delta-\delta 2^{-k-1}$
whose 
historical paths were displaced by distance more than
$\delta^{1/2-\eta}2^{-k/4}$
 on the time interval
$[a_n+\delta-\delta 2^{-(k-1)}, a_n+\delta
 -\delta 2^{-k}]$.
\begin{lemma}
\label{lem:nov2}
There exist $C=C(K_1)$
and $\delta_{\ref{lem:nov2}}$ such that,
for all $n$ sufficiently large,
\begin{equation}
\PP\left(\tZ^k_{a_n\,,\delta} >0\right) \leq
C
c_0e^{-\delta^{-\eta}2 ^{k/2}},\;\; \forall \delta\leq \delta_{\ref{lem:nov2}}.
\end{equation}
\end{lemma}
\paragraph{Proof:}
Let
$$\hZ^k_{a_n\,,\delta} = n\int_0^{\tau^{n,0}_{c_0}}
\mbox{\bf 1}_{ {\cal W}_{n,a,\delta,k,\eta,s}}
%
\ell^{n,a_n+\delta-\delta 2^{-k}}(ds),$$
that is,
$\hZ^{k}_{a_n\,,\delta}$ is  the total number of particles that are
alive at time $a_n+\delta-\delta 2^{-k}$ and whose historical paths
were
displaced by distance more than $\delta^{1/2-\eta}2^{-k/4}$
on the time interval
 $[a_n+\delta-\delta 2^{-(k-1)}, a_n+\delta
 -\delta 2^{-k}]$. We enumerate these particles by
 $i=1,\ldots,  \hZ^{k}_{a_n\,,\delta}$ and let
  $\hZ^{i,k}_{a_n\,,\delta}$ be the number of
  living descendents  of the particle $i$ ($i=1,\ldots,
 \hZ^{k}_{a_n\,,\delta}$) at time $a_n+\delta-\delta 2^{-k-1}$.
 Then clearly
 \begin{equation}
     \label{eq-nov26c}
     \tZ^k_{a_n\,,\delta} =
 \sum_{i=1}^{\hZ^{k}_{a_n\,,\delta}} \hZ^{i,k}_{a_n\,,\delta}.
\end{equation}
Lemma~\ref{lem:4} and (\ref{eq-nov26c}) imply that
for all $n$ sufficiently large
\begin{eqnarray*}
{\PP\left(\left.\tZ^k_{a_n\,,\delta} >0 \right|
X^{n,c_0}_{0,a_n+\delta(1-2^{-k})}\right)}
&\leq&
 \sum_{i=1}^{\hZ^{k}_{a_n\,,\delta}} \PP\left(\left.\hZ^{i,k}_{a_n\,,\delta} >0 \right|  X^{n,c_0}_{0,a_n+\delta(1-2^{-k})}\right)
\\
\nonumber
&\leq& \hZ^{k}_{a_n\,,\delta} 2h(b,\delta 2^{-k-1})/n
\leq \frac{4\hZ^{k}_{a_n\,,\delta}}{\delta 2^{-k-1} n},
\end{eqnarray*}
where the last inequality follows, for all $\delta$ sufficiently small,
from the definition of $h$.
Therefore,
$$\PP\left(\tZ^k_{a_n\,,\delta} >0\right)=
\EE\left( \PP\left(\left.\tZ^k_{a_n\,,\delta} >0 \right| X^{n,c_0}_{0,a_n+
\delta(1-2^{-k})}\right)\right)
\leq \frac{4}{\delta 2^{-k-1} n} \EE\left( \hZ^{k}_{a_n\,,\delta}\right).
$$
We next
represent the measure $X^{n,c_0}_{0,a_n+\delta(1-2^{-(k-1)})}$ as
\begin{eqnarray}
\label{equt:nov3}
 X^{n,c_0}_{0,a_n+\delta(1-2^{-(k-1)})}=\frac{1}{n} \sum_{i=1}^{nX^{n,c_0}_{0,a_n+\delta(1-2^{-(k-1)})}(1)}
  \delta_{\cz_i}
\end{eqnarray}
where $\cz_i$ are the
positions of the particles alive at time
$a_n+\delta(1-2^{-(k-1)})$. For the rest of the proof
of the lemma we call the particle that is located at $\cz_i$ at time
$a_n+\delta(1-2^{-(k-1)})$ --- the $i$-th particle.
Let $\tilde X^{n,i}$ be the measure describing the
positions of the living
descendents of the $i$-th particle at time
$a_n+\delta(1-2^{-k})$ and similarly to~(\ref{equt:nov3}) we can write
 \begin{eqnarray}
\tilde X^{n,i}=\frac{1}{n} \sum_{i=1}^{n\tilde X^{n,i}(1)}
  \delta_{\cz_{i,k}}
\end{eqnarray}
where $\cz_{i,k}$ is the position of the $k$-th descendent of the
$i$-th particle at time $a_n+\delta(1-2^{-k})$.
Then we get that
\[ X^{n,c_0}_{0,a_n+\delta(1-2^{-k})}= \sum_{i=1}^{nX^{n,c_0}_{0,a_n+\delta(1-2^{-(k-1)})}(1)}
   \tilde X^{n,i}.\]
Define
\[ f_{z}(x)= \mbox{\bf 1}_{|x-z|> \delta^{1/2-\eta}2^{-k/4}}\,, \;\; x,z\in \RR^d.\]
Then,
$$\hZ^{k}_{a_n\,,\delta}=
\sum_{i=1}^{nX^{n,c_0}_{0,a_n+\delta(1-2^{-(k-1)})}(1)}
   \tilde X^{i}_{a_n+\delta(1-2^{-(k-1)}, a_n+\delta(1-2^{-k})}(f_{\cz_i}).
$$
Hence, using
Lemma~\ref{lem:nov1} in the first inequality,
there exists $\delta_{\ref{lem:nov2}}$ sufficiently small such that
\begin{eqnarray*}
\lefteqn{\EE\left(\left. \hZ^{k}_{a_n\,,\delta}\right|
X^{n,c_0}_{0,a_n+\delta(1-2^{-(k-1)})}  \right)}
\\
&\leq&
\sum_{i=1}^{nX^{n,c_0}_{0,a_n+\delta(1-2^{-(k-1)})}(1)}
 \left(1+ \frac{\nu}{n}+\frac{\Vert \bar g\Vert_{\infty}}{2n}\right)^{n
 \delta 2^{-k}+1}\frac{1}{n}
 \PP_{\cz_i}(|B_{\delta 2^{-k}}-\cz_i| > \delta^{1/2-\eta} 2^{-k/4})\\
&\leq& X^{n,c_0}_{0,a_n+\delta(1-2^{-(k-1)})}(1) \left(1+
\frac{\nu}{n}+\frac{\Vert \bar g\Vert_{\infty}}{2n}\right)^{n\delta 2^{-k}+1}
e^{-\delta^{-\eta}2 ^{k/2}},\;\;\forall \delta\leq \delta_{\ref{lem:nov2}}\,,
\end{eqnarray*}
where $\PP_x$ is the law of the standard Brownian motion starting at $x$.
By taking the expectation we conclude
that for all $n$ sufficiently large,
\begin{eqnarray*}
\EE\left(\hZ^{k}_{a_n\,,\delta}\right)
&\leq&
C c_0 e^{-\delta^{-\eta}2 ^{k/2}}, \;\; \forall a\leq K_1, \delta\leq \delta_{\ref{lem:nov2}},
\end{eqnarray*}
where $C=C(K_1)$,
and we are done.
\gdm

\paragraph{Proof of Lemma~\ref{lem:nov4}:}
Fix $\delta_0$ sufficiently small such that $10^3\delta_0^{\eta/2}\leq 1$.
Let
$\delta\leq \delta_0$.
If $\tZ^k_{a_n\,,\delta}=0$ for each $k\geq 1$ then the maximal displacement
of the path of any particle on the time interval $[a_n, a_n+\delta]$
is bounded by
\[ \sum_{k=1}^{\infty} \delta^{1/2-\eta} 2^{-k/4} \leq
\delta^{1/2-\eta} \frac{1}{2^{1/4}-1}
 \leq  \delta^{1/2-1.5\eta}.\]
Hence by Lemma~\ref{lem:nov2} we get that for
$\delta\leq (\delta_0 \wedge \delta_{\ref{lem:nov2}})$,
$$\PP(Z^{n,3\eta/2}_{a_n\,,\delta}>0)\leq
\sum_{k=1}^{\infty}\PP(\tZ^k_{a_n\,,\delta}>0)
\leq C c_0  \sum_{k=1}^{\infty}e^{-\delta^{-\eta}2 ^{k/2}}
. $$
Now take $\delta_{\ref{lem:nov4}}\leq  (\delta_0 \wedge \delta_{\ref{lem:nov2}})$ sufficiently small so that for any $\delta\leq \delta_{\ref{lem:nov4}}$
\[  C
c_0 \sum_{k=1}^{\infty} e^{-\delta^{-\eta}2 ^{k/2}}\leq e^{-\delta^{-\eta}},\]
and we are done.
\gdm

\begin{lemma}
\label{lem:20_1} For any $\ep>0$, there exists 
$\delta_{1}>0$ such that
\begin{eqnarray}
\limsup_{n\rightarrow \infty} \PP\left(\sup_{a\leq K_1}\sup_{\delta\leq \delta_{1}} Z^{n,2\eta}_{a,\delta}>0\right) \leq \ep.
\end{eqnarray}
\end{lemma}

\paragraph{Proof:}
For any $m_0>0$ we have
by Lemma~\ref{lem:nov4} that
\begin{eqnarray*}
A_{m_0}&:=&\PP\left(Z^{n,3\eta/2}_{i2^{-m}, 2^{-m}}>0, {\rm for \;some}\;
i\leq K_12^m,\; m\geq m_0\right)\\
&\leq& \sum_{m=m_0}^{\infty} \sum _{i=0}^{K_12^m}
\PP\left(Z^{n,3\eta/2}_{i2^{-m}, 2^{-m}}>0\right) \leq
\sum_{m=m_0}^{\infty}K_12^m e^{-2^{m\eta}} \,.\end{eqnarray*} Choose
$m_0$ large enough so that $2^{-m_0}\leq \delta_{\ref{lem:nov4}}$,
$A_{m_0}\leq e^{-2^{m_0\eta/2}}\leq \ep$, and 
\begin{eqnarray}
\label{20_1}
10\cdot (2\cdot 2^{-m_0})^{1/2 -3\eta/2}
\leq \(2^{-m_0}\)^{1/2 -2\eta}.
\end{eqnarray}
Define
\[ C(K_1,m_0)=\{\omega:\; Z^{n,3\eta/2}_{i2^{-m}, 2^{-m}}=0,\;\;
\forall m>m_0, i\leq K_12^m\}.\]
Then
$$\PP(C(K_1,m_0))\geq 1-e^{-2^{m_0\eta/2}}
\geq 1-\ep.$$
Fix $\omega\in C(K_1,m_0)$. Fix arbitrary $a\leq K_1$ and
$\delta\leq 2^{-m_0}$. Then there exists $m\geq m_0$ such that
\begin{eqnarray}
\label{equt:nov5}
2^{-m-1} \leq \delta \leq 2^{-m}.
\end{eqnarray}
For $j\geq m_0$
let $\tilde a_j$ denote the smallest integer
multiple of $2^{-j}$ that is larger than
$a$ and, with $b=a+\delta$, let $\tilde b_j$ denote the largest integer
multiple
of $2^{-j}$ that is smaller than $b$.
Let $s$ be any time such that $Y^n_s =a+\delta$. Then
since
 $\delta\leq 2^{-m_0}\leq \delta_{\ref{lem:nov4}}$ and
$\omega\in C(K_1,m_0)$, we have by~(\ref{equt:nov5})
and the continuity of $\bW^n_s(\cdot)$
that
\begin{eqnarray*}
\left| \hat\WW^n_s - \bW^n_s(a)\right| &\leq& \left| \bW^n_s(\tilde
b_m) - \bW^n_s(\tilde a_m)\right|
+ \sum_{l=m+1}\left| \bW^n_s(\tilde a_l) - \bW^n_s(\tilde a_{l-1})\right|
\\
&&\mbox{}
+ \sum_{l=m+1}\left| \bW^n_s(\tilde b_l) - \bW^n_s(\tilde b_{l-1})\right|
\\
&\leq& 10\cdot 2^{-(1/2 -3\eta/2)m}
\leq 10\cdot (2\delta)^{1/2 -3\eta/2}
\leq \delta^{1/2 -2\eta},
\end{eqnarray*}
where the last inequality holds
by~(\ref{20_1}). By setting $\delta_{1}=2^{-m_0}$ we are done.
\gdm

The following corollary is immediate.
\begin{cor}
    \label{cor-immediate}
For any  $\ep>0$ there exists $\delta_1>0$ such that, 
\[ \PP\left(\sup_{s\leq \tau^{n,0}_{c_0}}\sup_{\delta\leq \delta_1}\sup_{a\leq (Y^n_s-\delta)_+}\left|
 \bW^n_s(a+\delta)-\bW^n_s(a)\right| > \delta_1^{1/2-2\eta}\right)\leq \ep.
\]
\end{cor}

We have made all the preparation for the proof of the following lemma,
concerning the tightness of the sequence $\{\WW^n\}_{n\geq 1}$.
\begin{lemma}
\label{lem:07_06_1} The sequence of processes $\{\WW^n\}_{n\geq 1}$ is $C$-tight in $D_\cW$.
\end{lemma} \paragraph{Proof:}
Recall that the $C$-tightness of the sequence of the contour
processes $\{Y^n\}_{n\geq 1}$ was proved in Section \ref{sec-tcont}
(see 
Proposition~\ref{prop:tight}).
Fix arbitrary $\beta>0$ and $\alpha=\beta^{1/2-2\eta}$. Then
for any $\delta_1>0$, we have the following inclusion 
\begin{eqnarray*}
    &&\{\sup_{s\leq \tau^{n,0}_{c_0}}\sup_{\delta\leq \delta_1}\sup_{u\geq 0} |\bW^n_{s+\delta}(u)-
    \bW^n_s(u)|\geq \alpha\}\subset\\
    &&
    \{\sup_{s\leq \tau^{n,0}_{c_0}}\sup_{\delta\leq \delta_1} |Y^n_{s+\delta}-Y^n_s|\geq \beta\}
    \\
    && \quad \bigcup
    \{\sup_{s\leq \tau^{n,0}_{c_0}} \sup_{\delta\leq \beta}
    \sup_{a\leq (Y_s-\delta)_+}
    |\bW^n_s(a+\delta)-\bW^n_s(a)|\geq \beta^{1/2-2\eta}\}.
\end{eqnarray*}
The $C$-tightness of the sequence $\{\WW^n\}_{n\geq 1}$ now follows
from this inclusion together with 
Proposition~\ref{prop:tight}, Corollary
\ref{cor-immediate}, and Lemma~\ref{lem:01_06_2}{\bf (b)}.
\gdm

We next turn to the local time processes $\ell^n, n\geq 1$.
\begin{lemma}
    \label{lem-newtl}
The sequence of processes  $\{\ell^{n,\cdot}\}_{n\geq 1}$ is
$C$-tight in $D_{\cM(\RR_+)}\,.$ 
\end{lemma}
\paragraph{Proof:}
Fix an arbitrary $c_0>0$, and define 
\[ \tilde \ell^{n,s}_t \equiv \ell^{n,s}_{t\wedge \tau^{n,0}_{c_0}},\; s,t\geq 0, \]
with $ \tilde \ell^{n,s}(dt)$ being as usual the corresponding
measure. 
Note that since
$c_0$ is arbitrary, it is enough to show the $C$-tightness  of
 $\{\tilde \ell^{n,\cdot}\}_{n\geq 1}$ in $D_{\cM_F}[0,\infty)$ and then the result follows immediately from
Lemma~\ref{lem:01_06_2}(b) (recall the
properties of convergence in vague topology). Since for each $t,n$, the function  $s\mapsto \tilde \ell^{n,t}_s$ is non-decreasing,
  to show the $C$-tightness  of
 $\{\tilde \ell^{n,\cdot}\}_{n\geq 1}$ in $D_{\cM_F}[0,\infty)$,
it is sufficient to prove the tightness of $\{\tilde \ell^{n,\cdot}_t\}_{n\geq 1}$
for each fixed $t$. That is, in view of Lemma
\ref{lem:01_06_2},
we need to prove that
for any constant $C$,
\begin{equation}
    \label{eq-technion0}
\limsup_{h\to 0}\limsup_{n\to\infty}
\PP(\sup_{0\leq r\leq C} |\tilde \ell^{n,r+h}_t-
\tilde \ell^{n,r}_t|>\epsilon)=0\,.
\end{equation}
The proof requires some care since introducing the time $t$ prevents
one from directly exploiting martingale properties and the tightness results
in \cite{M}.

We 
use the inverse local times
$\tau^{n,a}_r, a\geq 0, r\geq 0$
to define the
collection of processes
$$ \bar{X}_s^{i,j,\delta}=
\tilde \ell^{n,i\delta+s}_{\tau^{n,i\delta}_{(j+1)\delta}}-
\tilde \ell^{n,i\delta}_{\tau^{n,i\delta}_{j\delta}}\,,\;\;s\geq 0.$$
Note that $\bar{X}_s^{i,j,\delta}$ represents the total mass of the branching process in
random environment $X_{i\delta,i\delta+s}^{n,j\delta,(j+1)\delta}$, defined by~(\ref{eq:31_12_01}), which starts at
``time'' $i\delta$, such that
$$\bar{X}_0^{i,j,\delta}= X_{i\delta,i\delta}^{n,j\delta,(j+1)\delta}(1)=\delta.$$
 We also denote by ${\cal F}^{i,j,\delta}_l$ the filtration generated
by the process  $X_{i\delta,i\delta+\cdot}^{n,j\delta,(j+1)\delta}$ and its environment by time $l/n$.

On the event $t<\tau^{n,0}_{c_0}$ we have, for any $T>0$,
\begin{eqnarray}
    &&\sup_{0\leq r\leq T} |\tilde \ell^{n,r+h}_t-
\tilde \ell_t^{n,r}|\nonumber\\
&& \leq
\sup_{i\delta \leq T,j\delta \leq c_0}
\sup_{v\in [0, \delta]}
|\tilde \ell_{\tau^{n,i\delta}_{j\delta}}^{i\delta+v+h}
-\tilde \ell_{\tau^{n,i\delta}_{j\delta}}^{i\delta+v}|
+
\sup_{i\delta \leq T,j\delta \leq c_0}
\sup_{s\leq \delta} \bar{X}_s^{i,j,\delta}\nonumber\\
&&=:
\sup_{i\delta \leq T,j\delta \leq c_0}
A_{i,j}+
\sup_{i\delta \leq T,j\delta \leq c_0}
B_{i,j}\,.
\label{eq-technion1}
\end{eqnarray}
By the $C$-tightness of the sequence $\{s\mapsto \tilde
\ell_{\tau^{n,i\delta}_{j\delta}}^{i\delta+s}\}_{n\geq 1}$, see e.g. \cite{M}, Theorem
4.2 (proved there for the binary branching but valid, with similar proof,
for the geometric case under consideration here), we have that for each fixed
$\delta$ and each fixed $i\leq T/\delta, j\leq c_0/\delta$,
$$  \lim_{h\to 0}\limsup_{n\to\infty} \PP(A_{i,j}>\epsilon)=0\,.
$$
In particular, for any $\delta>0$ fixed,
\begin{equation}
\label{eq-technion5}
\lim_{h\to 0}\limsup_{n\to\infty}
    \PP(\sup_{i\delta \leq T,j\delta \leq c_0}
A_{i,j}>\epsilon)=0\,.
\end{equation}
To control $B_{i,j}$, we use the following lemma.
\begin{lemma}
    \label{lem-neqtech}
    For some universal constant $c$ and all $n$ large,
    $$
    \EE \sup_{0\leq s\leq \delta} ( (\bar{X}_s^{i,j,\delta})^4)\leq c \delta^4\,,
\quad \mbox{\rm for all $\delta\leq 1$}$$\end{lemma}
Indeed,  Lemma \ref{lem-neqtech} and Chebychev's inequality imply that
$$
\PP(\sup_{i\delta \leq T,j\delta \leq c_0} B_{i,j}>\epsilon)\leq
Tc_0 \delta^{-2}\delta^4\,.$$
Together with (\ref{eq-technion5}), this yields the proof of
Lemma \ref{lem-newtl},
once we complete the proof of Lemma \ref{lem-neqtech}.
\qed

In the proof of Lemma~\ref{lem-neqtech} we will frequently use the
following lemma, whose immediate proof (using iterations) is omitted.

\begin{lemma}
\label{recurs}
Let $c_1\,, c_2>0$ and suppose $z_i\,, i=1,2,\ldots$ satisfies the following inequalities
\[ z_i \leq \frac{c_1}{n}+ (1+\frac{c_2}{n}) z_{i-1}\,, \; i=1,2,\ldots.\]
Then  there exists $\bar{c}>0$ such that for any $\delta\in [0,1]$
\[ z_i\leq \bar{c}(\frac{c_1}{c_2}\delta + z_0),\;\; \forall i\leq \lfloor n\delta\rfloor.\]
\end{lemma}

\paragraph{Proof of Lemma \ref{lem-neqtech}:}
The argument uses computations similar
to those in  Section \ref{sec-asymp}. Throughout the proof, $\bar c$
denotes a constant whose value may change from line to line, but is independent
of $n$ or $\delta$.
Note that
the estimates on $\bar{X}_s^{i,j,\delta}$ that we get throughout the
proof below are uniform in $i,j$ and thus we may and will just
consider $i=j=1$ and write $\bar{X}_s=\bar{X}_s^{1,1,\delta}$ and
${\cal F}_l={\cal F}^{1,1,\delta}_l\,,l=0,1,2,\ldots$. Note that
$\bar{X}_s$ is the local time at level $s$ accumulated by the random
walk during its first 
$\lfloor n\delta\rfloor$ excursions from $0$.
We have the representation
$$ \bar{X}_{(m+1)/n}=n^{-1}\sum_{k=1}^{n\bar{X}_{m/n}} Z_{k,m+1}\,,$$
where the $Z_{k,m+1}$ is the number of offspring of the $k$-th particle at time $(m+1)/n$.
Recall that $Z_{k,m}\,, k=1,2,\ldots,$ are conditionally indpendent given ${\cal F}_m$, and for
each $k$, $Z_{k,m}$ is geometrically distributed with parameter
$1/2-\xi_{k,m}/4\sqrt{n}$. Here with some abuse of notation,
$$ \xi_{k,m}= \xi_{m/n}(x_{k,m}),$$
$\xi$ is as in
Section~\ref{1.1}, and $x_{k,m}$ is
the position of $k$-th particle at time $m$.
Note that by~(\ref{eqm32}) and our moment assumptions on $\xi$ we have that
$$\alpha_{k,m+1}:=\EE(Z_{k,m+1}|{\cal F}_m)\leq 1+\bar c/n.$$
Because the mean of $Z_{k,m}$ is close to $1$, the sequence
$\bar{X}_{(i+1)/n}$ is almost a martingale. To make it into a
martingale, introduce the variables, $M_0=\delta$,
$$ M_i=\frac{M_{i-1}}{\bar{X}_{(i-1)/n}}
\frac1n \sum_{k=1}^{n\bar{X}_{(i-1)/n}} \frac{Z_{k,i}}{\alpha_{k,i}},\;\;i\geq 1.$$
Note that
\begin{eqnarray}
\label{eqm33}
\bar{X}_{i/n}/M_i\leq (1+\bar c/n)^i,\;i\geq 1.
\end{eqnarray}
On the other hand, $i\mapsto M_i$ is a discrete martingale, and
hence by the Doob-Burkholder-Gundy
inequality, we have that
\begin{equation}
    \label{eq-wis1}
    \EE(\sup_{0\leq i\leq \delta n}
     M_i^4)\leq \bar c \EE\langle M\rangle_{\delta n}^2=
    \EE(\sum_{i=1}^{n\delta} \langle \Delta M\rangle_i)^2\,,
\end{equation}
where
$$\langle \Delta M\rangle_i=
\EE( (M_i-M_{i-1})^2|{\cal F}_{i-1})\,.$$

We prepare next some estimates.
First recall~(\ref{eqm32}), our moment assumptions on $\xi$
and its covariance structure
to get the following bound on the  correlation
between the $\{Z_{k,i+1}\}$:
    $$|\EE[(Z_{k,i+1}/\alpha_{k,i+1}-1)(Z_{k',i+1}/\alpha_{k',i+1}-1)|{\cal F}_{i}]|
    \leq \bar c/n\,,\;\;\forall k\not= k'.$$
Then we easily get,
\begin{eqnarray}
    \label{eq-wis11}
&&  \langle \Delta M\rangle_{i+1}=
    M_{i}^2 \EE\left( \left(\frac{1}{n\bar{X}_{i/n}}\sum_{k=1}^{n\bar{X}_{i/n}}
    \left(\frac{Z_{k,i+1}}{\alpha_{k,i+1}}-1\right)\right)^2|{\cal F}_{i}\right)\nonumber\\
    &\leq &
    \bar cM_i^2\frac{1}{n\bar{X}_{i/n}}+
    M_i^2 \max_{k\neq k', k,k' \leq n\bar X_{i/n}}
    E[(Z_{k,i+1}/\alpha_{k,i+1}-1)(Z_{k',i+1}/\alpha_{k',i+1}-1)|{\cal F}_{i}]\nonumber\\
    &\leq &
    \bar c\frac{M_i}{n}+\bar c \frac {M_i^2}{n}\,,
\end{eqnarray}
    Note that  $\EE M_i=\EE 
M_0=\delta$, and hence to control the right 
side of~(\ref{eq-wis11}) we need
to bound $\EE(M_i^2)$. $M_i$ is a martingale and hence  with $B_{1,i} = \EE(M_i^2)$ we use
~(\ref{eq-wis11}) to  get
 $$B_{1,i} \leq B_{i-1}+ \bar c\EE(\frac{M_{i-1}}{n})+
 \bar c \EE(\frac {M_{i-1}^2}{n})
 \leq (1+\frac{\bar c}{n}) B_{1,i-1} +\frac{\bar c \delta}{n}.
$$
By Lemma~\ref{recurs} we get
$$\EE(M_i^2)= B_{1,j}\leq \bar c( \delta^2 + M_{0}^2)
\leq  \bar c \delta^2,\;\; i\leq \lfloor n\delta\rfloor.
$$
Now recall again that  $\EE
M_i=\EE
M_0=\delta$ and use the above and~(\ref{eq-wis11}) to obtain that
    $$ \EE\langle M\rangle_i\leq \bar c\delta^2,\;\; \,i\leq \lfloor n\delta\rfloor.$$
    A similar computation, using Remark \ref{rem-new}, gives
        $$\EE( (M_{i+1}-M_i)^3| {\cal F}_i)
        \leq {\bar c}{n^{-2}} M_i+
        \bar c n^{-3/2} M_i^2+ \bar c n^{-1}M_i^3.$$
        With $B_{2,j}=E M_j^3$ one then obtains the recursions
        \begin{eqnarray*}
            B_{2,j+1}&\leq&\EE(M_j^3)+ \EE((M_{i+1}-M_i)^3)+ \bar c \EE(\EE(M_{i+1}-M_i)^2|{\cal F}_i)M_i)
\\
&\leq&
 (1+\frac{\bar c}{n}) B_{2,j}
        +\EE(M_j^2)(\bar c n^{-3/2} + \bar c n^{-1})
\nonumber
+ \EE(M_j) n^{-2} \\
\nonumber
&\leq&(1+\frac{\bar c}{n}) B_{2,j}
        +\bar c\delta^2 n^{-1}
\,,
\end{eqnarray*}
for $n$ sufficiently large ($n\geq \delta^{-1}$),
        and therefore by Lemma~\ref{recurs} we have
        \begin{equation}
            \label{eq-wis2}
            B_{2,j}\leq  \bar c( \delta^3 + M_{0}^3) \leq \bar c\delta^3, \;\; i\leq \lfloor n\delta\rfloor.
        \end{equation}
        Repeating this computation for the fourth moment, one
        obtains that
with $B_{3,j}=EM_j^4$,
        \begin{equation}
            \label{eq-wis3}
            B_{3,j}\leq  \bar c\delta^4, \;\; i\leq \lfloor n\delta\rfloor,
        \end{equation}
for all $n$ sufficiently large.
        Substituting (\ref{eq-wis11}) into (\ref{eq-wis1}) and using
        the last estimates, one gets
\begin{equation}
    \label{eq-wis10}
    \EE(\sup_{0\leq i\leq \delta n}
    M_i^4)\leq \bar c\delta^4\,,
\end{equation}
for all $n$ sufficiently large.
Since, by~(\ref{eqm33}),
$$
\sup_{0\leq s\leq \delta}  \bar{X}_s^4\leq \left( 1+\frac{\bar c}{n}\right)^{\delta n}
    \sup_{0\leq i\leq \delta n}
    M_i^4\,,
$$
this completes the proof of Lemma \ref{lem-neqtech}.
\qed

\begin{cor}
\label{cor:4}
$\{(\WW^n,\ell^n)\}_{n\geq 1}$ is $C$-tight in $D_{\cW\times\cM(\RR_+)}$.
\end{cor}
\paragraph{Proof:}
Immediately from Lemma \ref{lem-newtl}
and Lemma~\ref{lem:07_06_1}.
\qed

In what follows let
$(\bW,Y,\ell,\tau^0)$ be a limiting point of
$\{(\bW^n,Y^n,\ell,\tau^{n,0})\}_{n\geq 1}$. To simplify the
notation we 
omit subsequences
and simply
assume that $\{(\bW^n,Y^n,\ell^n,\tau^{n,0})\}_{n\geq 1}$
converges to $(\bW,Y,\ell,\tau^0)$.
We also switch (by Skorohod's theorem) to some probability space where
the convergence holds a.s..
Recall again that we write
$\ell^n_t$ and $\ell_t$ for
$\ell^n([0,t])$ and $\ell([0,t])$
respectively.

\begin{lemma}
    \label{lem-loctime}
$\ell$ 
is the
local time of $Y$.
\end{lemma}
\paragraph{Proof:}
First note that by properties of weak convergence of measures, for any $a\geq 0$
\begin{equation}
\label{01_06_3}
 \ell^{n,a}_t\rightarrow \ell^a_t
\end{equation}
for any point of continuity of function $t\mapsto \ell^a_t$. However
by a
limiting argument and the
convergence of $Y^n$ to $Y$, it is easy to
derive that if $Y_s\not= a$ then $s$ is a point of continuity of
$t\mapsto \ell^a_t$. Therefore, for all $a,t$ such that $Y_t\not= a$,
(\ref{01_06_3}) follows. Note that
$$
T^n_t(a)=
\frac{1}{n^2}
\sum^{\lfloor n^2 t\rfloor}_{i=0}{\bf 1}_{Y^n_{n^{-2}i}\le
a}
=
\int_0^{\lfloor n^2 t\rfloor/n^2}{\bf 1}_{Y^n_{s}\le a}
 \,ds, \;\; t\geq 0.
$$
Also  for any $a\geq 0$ and $\delta>0$  we have
$$
\int_{0}^{t} {\bf 1}_{a-\delta \leq Y^n_s\leq a+\delta}\,ds
=\int_{a-\delta}^{a+\delta} \ell^{n,s}_t\,ds
\leq 2\delta \sup_{s\leq K_1}  \ell^{n,s}_t.
$$
Since $\{\ell^{n,s}_t\}_{n\geq 1}$ is tight and $\delta$ was arbitrary we can
make the left side
arbitrarily small by taking $\delta>0$ sufficiently
small with probability as close to $1$ as we wish uniformly in $n$.
This, by a
standard argument, that also uses the convergence of $\{Y^n\}_{n\geq 1}$,
implies that
\begin{eqnarray}
\label{01_06_4}
\int_0^{\lfloor n^2 t\rfloor/n^2}{\bf 1}_{Y^n_{s}\le a}
\,ds \rightarrow  \int_0^{t}{\bf 1}_{Y_{s}\le a} \,ds
\end{eqnarray}
for any 
$a\ge 0, t\geq 0$. On the other
hand
$$T^n_t(a)= \int_0^a  \ell^{n,r}_t\,dr
\rightarrow  \int_0^a  \ell^{r}_t\,dr, t\geq 0,$$
where the last convergence follows by convergence of $\ell^{n,r}_t$ at all the points $r,t$ such that $Y_t\not= r$
(there is just one level $r$ such that $Y_t=r$). This and~(\ref{01_06_4}) yield
\begin{eqnarray}
    \int_0^{t}{\bf 1}_{Y_{s}\le a} \,ds = \int_0^a  \ell^{r}_t\,dr,\;\;  t\geq 0,
\end{eqnarray}
for all $a,r$, and hence $\ell^r_t$ is indeed the local time of $Y$,  for any  $t\geq 0$.
\qed

\begin{rem}
The above lemma and Corollary~\ref{cor:4} finish the proof of Proposition~\ref{prop:3}.
\end{rem}

The 
next two lemmas are essential for the proof of the
``charaterization of the limit points" part of Theorem~\ref{thm:1}.
First we prove the continuity of the local time at the level zero.
\begin{lemma}
\label{lem:29_12_2} $t\mapsto \ell^0_t$ is continuous.
\end{lemma}
\paragraph{Proof:}
It is enough to show that for arbitrary $c_0>0$, $\left\{\ell^{n,0}_{\cdot\wedge \tau^{n,0}_{c_0}}\right\}_{n\geq 1}$ is $C$-tight in $D_{\RR}[0,\infty)$,
that is, for any $\ep>0$
\begin{eqnarray}
\label{eq:27_1}
\lim_{\delta\downarrow 0} \limsup_{n\rightarrow \infty}\PP\( \sup_{t\leq \tau^{n,0}_{c_0}}
\ell^{n,0}_t-\ell^{n,0}_{t-\delta}\geq \ep\)=0.
\end{eqnarray}
Suppose~(\ref{eq:27_1}) does not hold, that is, there exist $\ep, \ep_1>0$, such that for all
 $\delta>0$
\begin{eqnarray}
\label{eq:27_2}
\PP\( \sup_{t\leq \tau^{n,0}_{c_0}}
\ell^{n,0}_t-\ell^{n,0}_{t-\delta}\geq \ep\)\geq \ep_1.
\end{eqnarray}
Fix such $\ep, \ep_1>0$;
we have the inclusion
\[ \left\{ \sup_{t\leq \tau^{n,0}_{c_0}} \ell^{n,0}_t-\ell^{n,0}_{t-\delta}\geq \ep\right\}
 \subset \left\{\exists i=1,\ldots,\left\lfloor
\frac{2c_0}{\ep}\right\rfloor :\;  \tau^{n,0}_{\frac{(i+1)\ep}{2}} -
 \tau^{n,0}_{\frac{i\ep}{2}}<\delta\right\}. \]
Since $\tau^{n,0}_{\frac{(i+1)\ep}{2}} -
 \tau^{n,0}_{\frac{i\ep}{2}}, i=1,\ldots, \left\lfloor \frac{2c_0}{\ep}\right\rfloor$ are identically distributed we get
$$\PP\( \exists i=1,\ldots,\left\lfloor
\frac{2c_0}{\ep}\right\rfloor  :\;  \tau^{n,0}_{\frac{(i+1)\ep}{2}} -
 \tau^{n,0}_{\frac{i\ep}{2}}<\delta\)
\leq \(\left\lfloor  \frac{2c_0}{\ep}\right\rfloor+1\)
 \PP\(  \tau^{n,0}_{\ep/2}<\delta\).
$$
By Lemma~\ref{lem:01_06_2}(a), we can choose $\delta$ sufficiently small such that
\[  \PP\(  \tau^{n,0}_{\ep/2}<\delta\)\leq \frac{\ep_1}{2\(\left\lfloor  \frac{2c_0}{\ep}\right\rfloor+1\)}
\]
for all $n$ sufficiently large, and hence
\[ \PP\(\sup_{t\leq \tau^{n,0}_{c_0}} \ell^{n,0}_t-\ell^{n,0}_{t-\delta}\geq \ep\)\leq \frac{\ep_1}{2}\]
and we get a contradiction with (\ref{eq:27_2}).
\qed

\begin{lemma}
\label{lem:28_12_1} For any fixed $r>0$,  $\tau^0_r$ equals, with
probability one, to the value of the inverse function of
$\ell^0_{\cdot}$ at $r$, that is, 
\[ \tau^0_r=\inf\{s>0:\;\ell^0_s>r\}\,, \quad a.s..\]
\end{lemma}
\paragraph{Proof:}
Recall that we assume that we are considering the probability space
where $\ell^{n,0},\tau^{n,0}\rightarrow (\ell^0,\tau^0)$ in
$D_{\RR_+}[0,\infty)\times \cM(\RR_+)$, $\PP$-a.s.. Moreover we know
that for any fixed $r$, $\tau^0(\cdot)$ is continuous at the point
$r$. This, by 
properties of
convergence in $\cM$, implies that for any fixed $r$,
$\tau^{n,0}_r\rightarrow \tau^0_r$, $\PP$-a.s.. Fix arbitrary
$c_0,\delta>0$. Then, by definition of the local time,
 we get,
\begin{eqnarray}
\label{eq:29_12_1}
  \ell^{n,0}_{\tau^{n,0}_{c_0+\delta}}\geq c_0+\delta.
\end{eqnarray}

Since   $\ell^{n,0}_{\cdot}$ converges to the continous limit, the convergence is uniform on the compacts. This and
the convergence $\tau^{n,0}_{c_0+\delta}\rightarrow \tau^0_{c_0+\delta}$ imply, that by passing to the limit in~(\ref{eq:29_12_1}) we get
\begin{eqnarray}
\label{eq:29_12_2}
  \ell^{0}_{\tau^0_{c_0+\delta}}\geq c_0+\delta,
\end{eqnarray}
and hence
\begin{eqnarray}
\inf\{s>0:\;\ell^0_s>c_0\}\leq  \tau^0_{c_0+\delta}.
\end{eqnarray}
Similarly we can show that
\begin{eqnarray}
\inf\{s>0:\;\ell^0_s>c_0\}\geq  \tau^0_{c_0-\delta}.
\end{eqnarray}
Since $\delta$ was arbitrary, and by the continuity of $\tau^0$ at $c_0$
 (see Lemma~\ref{lem:01_06_2}{\bf (d)}) we get
\begin{eqnarray}
\inf\{s>0:\;\ell^0_s>c_0\}= \tau^0_{c_0}.
\end{eqnarray}
and we are done.

\qed

\begin{lemma}
\label{lem:27_12_3}
For any $\phi\in \cC_b(\RR^d)$  and fixed $c_0>0$,
\begin{eqnarray}
\int_0^{\tau^{n,0}_{c_0}}\phi(\hat \WW^n_s)\ell^{n,t}(ds) \rightarrow  \int_0^{\tau^0_{c_0}}\phi(\hat \WW_s)\ell^{t}(ds),
\;\;\forall t\geq 0,\; \PP-{\rm a.s.},\;
\end{eqnarray}
{\rm as}$\;n\rightarrow\infty,$ where
\begin{eqnarray}
\label{eq:29_12_4}
\tau^0_{c_0}=\inf\{r>0:\;\ell^0_r>c_0\}.
\end{eqnarray}
\end{lemma}
\paragraph{Proof:}
 $\tau^{n,0}_{c_0}\rightarrow \tau^0_{c_0}$, where by Lemma~\ref{lem:28_12_1} $\tau^0_{c_0}$ is defined by~(\ref{eq:29_12_4}). Moreover,
 by Lemma~\ref{lem:29_12_2}, $\ell^0_{\cdot}$ is continuous
at $ \tau^0_{c_0}$, therefore by elementary properties of weak convergence, for any continnuous function $f(s)$
\[ \int_0^{\tau^{n,0}_{c_0}}f(s)\ell^{n,0}(ds) \rightarrow  \int_0^{\tau^0_{c_0}}f(s)\ell^0(ds),\; \PP-{\rm a.s.},\;
\;{\rm as}\;n\rightarrow\infty.
\]
Now the result for $t=0$, follows by uniform on the compacts
convergence of $\hat \WW^n$ to $\hat \WW$. The convergence of the
integral  for $t>0$ follows immediately since, 
by the continuity of $Y$, the $\ell^t(ds)$ does not charge
the point $s=\tau^0_{c_0}$ for every $t>0$. \qed

\paragraph{Proof of Theorem \ref{thm:1}:}
The tightness statement was
proved in Proposition \ref{prop:3}. To finish the proof we need to
derive the characterization of the limit points. Fix arbitrary
$c_0>0$ and let $ X^{n,c_0}_{0,t}$ be the measure-valued process
defined as in~(\ref{eq:Oc5_1}), that is,
\begin{eqnarray}
\label{20_3}
X^{n,c_0}_{0,t}(\phi)&\equiv&\int_{0}^{\tau^{n,0}_{c_0}}\phi(\hat\WW^n_s)
\ell^{n,t}(ds),\;\;\; t\in [0,K_1],
\end{eqnarray}
for all $\phi\in \cB(\RR^d)$. Let $(\bW, Y, \ell, \tau^0_{c_0})$ be an 
arbirary limit point of \\
$\{(\bW^n, Y^n, \ell^n, \tau^{n,0}_{c_0})\}_{n\geq 1}$.
Fix arbitrary $\phi \in \cC_b(\RR^d)$.
As we have mentioned already, due to results in~\cite{M}, the sequence of process  $\{X^{n,c_0}_{0,\cdot}\}_{n\geq 1}$ converges weakly in $D_{\cM_F}[0,K_1]$ to the process $X^{c_0}\in C_{\cM_F}[0,K_1]$ satisfying the
martingale problem (\ref{MP-1}-\ref{MP-2}) on $[0,K_1]$,
with
 $X^{c_0}_0=c_0\delta_x$, and hence the left hand side of~(\ref{20_3}) converges to $X^{c_0}_t(\phi)$ for any $t\in [0,K_1]$.
As for the right hand side of~(\ref{20_3}), due to
Proposition~\ref{prop:3} and Lemma~\ref{lem:27_12_3} it converges,
along an appropriate subsequence, to
$\int_{0}^{\tau^{0}_{c_0}}\phi(\hat \WW_s) \ell^{t}(ds)$
for  $t\in
[0,K_1]$, where $\ell$ is  the local time $Y$. This gives
us~(\ref{20_2}) for any $\phi \in \cC_b(\RR^d)$. The extension of
the equality to any $\phi\in \cB(\RR^d)$ is trivial.
\qed

\section{Proof of Theorem~\ref{thr:5.1}}
\label{sec-5}
\setcounter{equation}{0}

The proof of the result is based on convergence of approximations.
For simplicity, as before, we  assume that
$(\WW^n,B^n,\ell^n)=(\bW^n,Y^n,B^n,\ell^n)\rightarrow
(\bW,Y,B,\ell)=(\WW,Y,B,\ell)$ a.s. (based on
Proposition~\ref{prop:3} we can always get it by Skorohod's theorem
via an appropriate subsequence). Further, we localize the snake $\bW^n$ to 
live in a compact, and then it is not hard to check that
for $n$ large, the truncation in the definition of $\xi_j^n(y)$ can be ignored. 
Thus, we assume in the sequel that
\[\frac{\xi_j^n\(y\)}{\sqrt n}= 
\frac{\xi_j\(y\)}{\sqrt n}= B_{\frac{j}{n}}(y)-B_{\frac{j-1}{n}}(y).\]
On the level of the $n$th approximation we will be
dealing with the following approximating functional:

$$ F_n\(\WW^n_{n^{-2}k}\)\equiv\frac{1}{n}
\sum_{l=1}^{nY^n_{k/n^2}-1} e^{-\frac{1}{\sqrt{n}}
\sum_{l'=1}^{l}\xi_{l'}\(\bW^{n}_{k/n^2}(\frac{l+1}{n})\)}. $$


Note that
\begin{eqnarray*}
F_n\(\WW^n_{n^{-2}(k+1)}\) &=&\left\{\begin{array}{ll} \frac{1}{n}
\sum_{l=1}^{nY^n_{k/n^2}-2}
 e^{-\frac{1}{\sqrt{n}}\sum_{l'=1}^{l}\xi_{l'}\(\bW^{n}_{k/n^2}(\frac{l+1}{n})\)
},&\mbox{if }
Y^n_{(k+1)/n^2}< Y^n_{k/n^2},\\
&\\
\frac{1}{n} \sum_{l=1}^{nY^n_{k/n^2}-1} e^{-\frac{1}{\sqrt{n}}
\sum_{l'=1}^{l}\xi_{l'}\(\bW^{n}_{k/n^2}(\frac{l+1}{n})\)}&
\\
\mbox{}
+\frac{1}{n}e^{-\frac{1}{\sqrt{n}}
\sum_{l'=1}^{nY^n_{k/n^2}}\xi_{l'}\(\bW^{n}_{(k+1)/n^2}(Y^n_{k/n^2}+1/n)\)},
&\mbox{if } Y^n_{(k+1)/n^2}> Y^n_{k/n^2}.\end{array}\right.
\end{eqnarray*}
Further, 
if
$Y^n_{(k+1)/n^2}> Y^n_{k/n^2}$, then
$$\bW^{n}_{(k+1)/n^2}(Y^n_{k/n^2}+1/n)=\hat{\WW}^n_{(k+1)/n^2}=
\hat{\WW}^n_{k/n^2}+\eta_{1/n},$$ where  $\eta$ is a Brownian path
independent of $\WW^n_{k/n^2}$. Let
\begin{eqnarray*}
\cF_k&=& \sigma\left\{\WW^n_{l/n^{2}}\,,\; l\leq k\right\}\vee \sigma\left\{\xi_{l}\,,\; l=0,1,2,\ldots\right\}.
\end{eqnarray*}

Define
\begin{eqnarray*}
V_{k}=F\(\WW^n_{k/n^2}\)-F\(\WW^n_{(k-1)/n^2}\),\;\; k=1,2,\ldots.
\end{eqnarray*}
Then by the standard decomposition  of $F\(\WW^n_{m/n^2}\)$ we get that
\begin{eqnarray*}
F\(\WW_{m/n^2}\)=M^n_m+ A^n_m,\;\;m=1,2,\ldots,
\end{eqnarray*}
where $M^n_m\,, m=1,2,\ldots,$ is the $\{\cF_m\}_{m\geq
1}$-martingale 
given by
\begin{eqnarray*}
M^n_m&=& \sum_{k=0}^{m-1}\left( V_{k+1}- \EE\left(V_{k+1}|
\cF_k\right)\right),\;\;m=1,2,\ldots
\end{eqnarray*}
and
\begin{eqnarray*}
A^n_m&=& \sum_{k=0}^{m-1}\EE\left(V_{k+1}| \cF_k\right),\;\;m=1,2,\ldots.
\end{eqnarray*}

We first study the limiting behavior of $A^n$.
\begin{lemma}
\label{lem:5.1}
\begin{eqnarray}
\nonumber
A^n_{\lfloor n^2 t\rfloor} &\rightarrow&
 \int_0^t e^{-B_{Y_s}(\hat{\WW}_s)}\left\{-\frac{1}{2}\Delta B_{Y_s}(\hat{\WW}_s)+\frac{1}{2}\sum_{i=1}^d
  \left(\frac{\partial}{\partial x_i} B_{Y_s}(\hat{\WW}_s)\right)^2\right\}\,ds \\
\label{eq:5.2}
&&+ \ell^0_t -\int_0^t e^{-B_{K_1}(\hat{\WW}_s)}\,\ell^{K_1}(ds),\;\;{\rm as}\; n\rightarrow\infty.
\end{eqnarray}
\end{lemma}
\paragraph{Proof:}
Using $\EE_\eta$ to denote expectation with respect to the Brownian
path $\eta_\cdot$, we have 
\begin{eqnarray*}
\lefteqn{\EE\(V_{k+1}|\cF_k\)}\\
&=&\PP\(Y^n_{(k+1)/n^2}< Y^n_{k/n^2}|\cF_k\)\\
&&\;\;\;\mbox{}\times
\EE\left(-\frac{1}{n}e^{-\frac{1}{\sqrt{n}}
\sum_{l'=1}^{nY^n_{k/n^2}-1}\xi_{l'}(\hat{\WW}^n_{n^{-2}k})}
\bigg|Y^n_{(k+1)/n^2}< Y^n_{k/n^2}\,,\  \cF_k\right)\\
&&+\PP\(Y^n_{(k+1)/n^2}> Y^n_{k/n^2}|\cF_k\)\\
&&\;\;\;\mbox{}\times
\EE\left(\frac{1}{n}e^{-\frac{1}{\sqrt{n}}\sum_{l'=1}^{nY^n_{k/n^2}}\xi_{l'}
\(\bW^{n}_{n^{-2}(k+1)}\big[Y^n_{k/n^2}+1/n\big]\)}
\bigg|Y^n_{(k+1)/n^2}> Y^n_{k/n^2}\,,\ \cF_k\right)\,.
\end{eqnarray*}
Therefore,
\begin{eqnarray*}
\lefteqn{\EE\(V_{k+1}|\cF_k\)}\\
&=&-\left[\frac{1}{2}-\frac{1}{4\sqrt{n}}\xi_{Y^n_{k/n^2}}\left(
\hat{\WW}^n_{n^{-2}k}\right) \right]
\frac{1}{n}e^{-\frac{1}{\sqrt{n}}\sum_{l'=1}^{nY^n_{k/n^2}-1}\xi_{l'}(\hat{\WW}^{n}_{n^{-2}k})}\\
&&\mbox{}+\left[\frac{1}{2}+\frac{1}{4\sqrt{n}}\xi_{Y^n_{k/n^2}}\left(
\hat{\WW}^n_{n^{-2}k}\right) \right]
\EE_{\eta}\(\frac{1}{n}e^{-\frac{1}{\sqrt{n}}\sum_{l'=1}^{nY^n_{k/n^2}}\xi_{l'}(\hat{\WW}^{n}_{n^{-2}k}+\eta_{1/n})}\)\\
&&\mbox{}+\frac{1}{n}{\bf 1}_{\{Y^n_{k/n^2}=0\}}
-{\bf 1}_{\{
Y^n_{k/n^2}=K_1\}}\frac{1}{n}e^{-\frac{1}{\sqrt{n}}\sum_{l'=1}^{nK_1-1}\xi_{l'}(\hat{\WW}^{n}_{n^{-2}k})}
\\
&=&\frac{1}{2n}\left[\EE_{\eta}\(e^{-\frac{1}{\sqrt{n}}\sum_{l'=1}^{nY^n_{k/n^2}}\xi_{l'}(\hat{\WW}^{n}_{n^{-2}k}+\eta_{1/n})}\)-
e^{-\frac{1}{\sqrt{n}}\sum_{l'=1}^{nY^n_{k/n^2}-1}
\xi_{l'}(\hat{\WW}^{n}_{n^{-2}k})}\right]\\
&&\mbox{}+\frac{1}{4n^{3/2}}\xi_{Y^n_{k/n^2}}\left(
\hat\WW^n_{n^{-2}k}\right) \left[e^{-\frac{1}{\sqrt{n}}
\sum_{l'=1}^{nY^n_{k/n^2}-1}\xi_{l'}(\hat{\WW}^{n}_{n^{-2}k})}
\right.\\
&&\left. \hspace*{2cm}\mbox{}
+
\EE_{\eta}\(e^{-\frac{1}{\sqrt{n}}\sum_{l'=1}^{nY^n_{k/n^2}}\xi_{l'}(\hat{\WW}^{n}_{n^{-2}k}+\eta_{1/n})}\)
\right]
\\
&&\mbox{}+\frac{1}{n}\cdot {\bf  1}_{\{Y^n_{k/n^2}=0\}}
-\frac{1}{n}\cdot {\bf 1}_{\{Y^n_{k/n^2}=K_1\}}e^{-\frac{1}{\sqrt{n}}\sum_{l'=1}^{nK_1-1}\xi_{l'}(\hat{\WW}^{n}_{n^{-2}k})}
\\
&=&\frac{1}{2n}\left[\EE_{\eta}\(e^{-B_{Y^n_{k/n^2}}(\hat{\WW}^{n}_{n^{-2}k}+\eta_{1/n})}\)-
e^{-B_{Y^n_{k/n^2}-1/n}(\hat{\WW}^{n}_{n^{-2}k})}\right]\\
&&\mbox{}+\frac{1}{4n^{3/2}}\xi_{Y^n_{k/n^2}}\left(
\hat{\WW}_{n^{-2}k}\right) \left(e^{-B_{Y^n_{k/n^2}-1/n}(\hat{\WW}^{n}_{n^{-2}k})}+
\EE_{\eta}\(e^{-B_{Y^n_{k/n^2}}(\hat{\WW}^{n}_{n^{-2}k}+\eta_{1/n})}\)
\right)
\\
&&\mbox{}+\frac{1}{n}\cdot {\bf 1}_{\{Y^n_{k/n^2}=0\}}
-\frac{1}{n}\cdot {\bf 1}_{\{Y^n_{(k-1)/n^2}=K_1-1,Y^n_{k/n^2}=K_1\}
}e^{-B_{K_1-1/n}(\hat{\WW}^{n}_{n^{-2}k})}\\
&=&I_{1,n,k}+I_{2,n,k}+I_{3,n,k}-I_{4,k,n}\,,
\end{eqnarray*}
where we also used the definition of $B$.
We begin with an estimate of
$\EE_{\eta}\(e^{-B_{Y^n_{k/n^2}}(\hat{\WW}^{n}_{n^{-2}k}+\eta_{1/n})}\)$.
Set
\begin{eqnarray*}
R_{n,k}(s)&=& e^{-B_{Y^n_{k/n^2}}(\hat{\WW}^{n}_{n^{-2}k}+\eta_{s})}\left[
  -\frac{1}{2}\Delta_x B_{Y^n_{k/n^2}}(\hat{\WW}^{n}_{n^{-2}k}+x)|_{x=\eta_{s}}
\right.\\
&& \left.\mbox{}\quad\quad+
\frac{1}{2}\sum_{i=1}^d\left(\frac{\partial}{\partial x_i}
B_{Y^n_{k/n^2}}(\hat{\WW}^{n}_{n^{-2}k}+x)|_{x=\eta_{s}}\right)^2
\right]\,.
\end{eqnarray*}

By It\^{o}'s formula we get
$$
\EE_{\eta}\(e^{-B_{Y^n_{k/n^2}}(\hat{\WW}^{n}_{n^{-2}k}+\eta_{1/n})}\)=
e^{-B_{Y^n_{k/n^2}}(\hat{\WW}^{n}_{n^{-2}k})}
+\EE_{\eta}\left(\int_{0}^{1/n}R_{n,k}(s)ds\right)\,.
$$
The first term at the right side above can be further
decomposed as
\begin{eqnarray*}
\lefteqn{e^{-B_{Y^n_{k/n^2}}(\hat{\WW}^{n}_{n^{-2}k})}}\\
&=& e^{-B_{Y^n_{k/n^2}-1/n}(\hat{\WW}^{n}_{n^{-2}k})}\left[
1-\frac{1}{n^{1/2}}\xi_{Y^n_{k/n^2}}\left(
\hat{\WW}_{n^{-2}k}\right)+\frac{1}{2n}\xi_{Y^n_{k/n^2}}\left(
\hat{\WW}_{n^{-2}k}\right)^2\right.\\
&&\left. \mbox{} +O(n^{-3/2})\Theta\left(\left|\xi_{Y^n_{k/n^2}}\left(
\hat{\WW}_{n^{-2}k}\right)\right|^3\right)\right]\,,
\end{eqnarray*}
where $\Theta(x)$ is some point in $[-x,x]$.
We get
\begin{eqnarray*}
I_{1,n,k}&=& \frac{1}{2n}\left[
e^{-B_{Y^n_{k/n^2}-1/n}(\hat{\WW}^{n}_{n^{-2}k})}\left(
-\frac{1}{n^{1/2}}\xi_{Y^n_{k/n^2}}\left(
\hat{\WW}_{n^{-2}k}\right)+\frac{1}{2n}\xi_{Y^n_{k/n^2}}\left(
\hat{\WW}_{n^{-2}k}\right)^2\right.\right.\\
&&\left. \mbox{} +O(n^{-3/2})\Theta\left[\left|\xi_{Y^n_{k/n^2}}\left(
\hat{\WW}_{n^{-2}k}\right)\right|^3\right]\right)
\\
&&\mbox{}+
\left.\EE_{\eta}\left(\int_{0}^{1/n}R_{n,k}(s) ds\right)\right]
. \end{eqnarray*}
To handle
$I_{2,n,k}$, 
note that
\begin{eqnarray*}
I_{2,n,k}&=& \frac{1}{4n^{3/2}}\xi_{Y^n_{k/n^2}}\left(
\hat{\WW}_{n^{-2}k}\right)\left[
e^{-B_{Y^n_{k/n^2}-1/n}(\hat{\WW}^{n}_{n^{-2}k})}\left(
2-\frac{1}{n^{1/2}}\xi_{Y^n_{k/n^2}}\left(
\hat{\WW}_{n^{-2}k}\right)\right.\right.
\\
&&\left. \mbox{} +O(n^{-1})\Theta\left[\left|\xi_{Y^n_{k/n^2}}\left(
\hat{\WW}_{n^{-2}k}\right)\right|^2\right]\right)
+\left.\EE_{\eta}\left(\int_0^{1/n} R_{n,k}(s)\,ds\right)\right].
\end{eqnarray*}
All together we get
\begin{eqnarray*}
I_{1,n,k}+ I_{2,n,k}&=&
\frac{1}{2n}\EE_{\eta}\left(\int_0^{1/n}R_{n,k}(s)\,ds\right) +O(n^{-5/2})
e^{-B_{Y^n_{k/n^2}-1/n}(\hat{\WW}^{n}_{n^{-2}k})}
\Theta\left(\left|\xi_{Y^n_{k/n^2}}\left(
\hat{\WW}_{n^{-2}k}\right)\right|^3\right)
\\
&&
\!\!\!\!\!\!\!\!\!
\!\!\!\!\!\!\!\!\!
\!\!\!\!\!\!\!\!\!
\mbox{}
+O(n^{-3/2})\Theta\left(\left|\xi_{Y^n_{k/n^2}}\left(
\hat{\WW}_{n^{-2}k}\right)\right|\right)
\EE_{\eta}\left(\int_0^{1/n} R_{n,k}(s)\,ds\right)\,.
\end{eqnarray*}
From this it follows that for any $t>0$
\begin{eqnarray*}
\sum_{k=1}^{\lfloor n^2t\rfloor}(I_{1,n,k}+ I_{2,n,k})&=&
\frac{1}{2n^2}\sum_{k=1}^{\lfloor n^2t\rfloor} n\EE_{\eta}\left(\int_0^{1/n} R_{n,k}(s)\,ds\right)\\
&&\mbox{} +
 O(n^{-5/2})\sum_{k=1}^{\lfloor n^2t\rfloor} \Theta\left(\left|\xi_{Y^n_{k/n^2}}\left(
\hat{\WW}_{n^{-2}k}\right)\right|^3\right)
e^{-B_{Y^n_{k/n^2}-1/n}(\hat{\WW}^{n}_{n^{-2}k})}
\\
&&
\!\!\!\!\!\!\!\!\!
\!\!\!\!\!\!\!\!\!
\!\!\!\!\!\!\!\!\!
\!\!\!\!\!\!\!\!\!
\!\!\!\!\!\!\!\!\!
+O(n^{-5/2})\sum_{k=1}^{\lfloor n^2t\rfloor}\Theta\left(\left|\xi_{Y^n_{k/n^2}}\left(
\hat{\WW}_{n^{-2}k}\right)\right|\right)
n\EE_{\eta}\left(\int_0^{1/n} R_{n,k}(s)\,ds\right)\\
&&
\!\!\!\!\!\!\!\!\!
\!\!\!\!\!\!\!\!\!
\!\!\!\!\!\!\!\!\!
\rightarrow
\int_0^t e^{-B_{Y_s}(\hat{\WW}_{s})}\left(
  -\frac{1}{2}\Delta_x B_{Y_s}(\hat{\WW}_{s})
+\frac{1}{2}\sum_{i=1}^d\left(\frac{\partial}{\partial x_i}
B_{Y_s}(\hat{\WW}_{s})\right)^2
\right)\,ds,
\end{eqnarray*}
where, due to the H\"{o}lder continuity of the Gaussian field,
the second and third  terms on
the right side of the first equality
converge to $0$ (once $\hat \WW_{n^{-2}k}$ has been localized 
in a compact with high probability),
and the
first term converges to the first term on the right side of~(\ref{eq:5.2}).

Now we will treat $I_{3,n,k}$ and $I_{4,n,k}$. By definition of the
approximate local time $\ell^{n,n^{-1}m}$ we get 
\begin{eqnarray*}
\sum_{k=0}^{\lfloor n^2t\rfloor}(I_{3,n,k}+ I_{4,n,k})&=&
\ell^{n,0}_t
- \int_0^t \EE_{\eta}\left(e^{-B_{K_1-\frac{1}{n}}(\hat{\WW}^{n}_{s+1/n^2})}
\right) \ell^{n,K_1-1/n}(ds).
\end{eqnarray*}
Then pass to the limit, use the 
uniform on compacts convergence of $\ell^n$
and $\WW^n$ to $\ell$ and $\WW$, 
and the continuity of $B$
to get
that
\begin{eqnarray*}
\sum_{k=0}^{\lfloor n^2t\rfloor}(I_{3,n,k}+ I_{4,n,k})&\rightarrow &   \ell^{0}_t - \int_0^t e^{-B_{K_1}(\hat{\WW}_{s})}\ell^{K_1}(ds),
\end{eqnarray*}
as $n\rightarrow\infty$.
Thus,
we obtain the second and the third terms in (\ref{eq:5.2}).
\qed

Define the bracket process for the martingale $M^n$:
\begin{eqnarray}
\<M^n_\cdot\>_m&\equiv& \sum_{k=0}^{m-1}
\EE\left(\left(M^n_{k+1}-M^n_{k}\right)^2|\cF_k\right),\;\;
 m=1,2,\ldots.
\end{eqnarray}

Then we have
\begin{lemma}
\label{lem:5.2}
\begin{eqnarray*}
 \langle M^n_{\cdot}\rangle_{\lfloor n^2 t\rfloor} \rightarrow
 \int_0^t e^{-2B_{Y_s}(\hat{\WW}_{s})}ds,\;\;{\rm as}\;n\rightarrow \infty.
\end{eqnarray*}

\end{lemma}

\paragraph{Proof:}
It is easy to check that for any $m\geq 1$,
\begin{eqnarray*}
 \langle M^n_{\cdot}\rangle_{m} &=&
\sum_{k=0}^{m-1} \EE\left(V_{k+1}^2|\cF_k\right) -
\sum_{k=0}^{m-1} \left(\EE\left(V_{k+1}|\cF_k\right)\right)^2.
\end{eqnarray*}
By Lemma~\ref{lem:5.1} we know that as $n\rightarrow\infty$
\begin{eqnarray*}
\sum_{k=0}^{\lfloor n^2 t\rfloor} \EE\left(V_{k+1}|\cF_k\right)&\rightarrow&
 \int_0^t e^{-B_{Y_s}(\hat{\WW}_s)}\left\{-\frac{1}{2}\Delta B_{Y_s}(\hat{\WW}_s)+\frac{1}{2}\sum_{i=1}^d
  \left(\frac{\partial}{\partial x_i} B_{Y_s}(\hat{W}_s)\right)^2\right\}\,ds \\
\nonumber
&&+ \ell^0_t -\int_0^t e^{-B_{K_1}(\hat{\WW}_s)}\,\ell^{K_1}(ds)
\end{eqnarray*}
which is a
process of bounded variation. From this it is easy to deduce that
\begin{eqnarray*}
\sum_{k=0}^{\lfloor n^2 t\rfloor} \left(\EE\left(V_{k+1}|\cF_k\right)\right)^2&\rightarrow&0,
\end{eqnarray*}
as $n\rightarrow \infty$.
Hence it is enough to consider the limiting behavior of
\begin{eqnarray*}
\sum_{k=0}^{\lfloor n^2 t\rfloor} \EE\left(V_{k+1}^2|\cF_k\right).
\end{eqnarray*}
By repeating the argument in the proof of Lemma~\ref{lem:5.1} we get
\begin{eqnarray*}
\lefteqn{\EE(V_{k+1}^2|\cF_k)}\\
&=&\PP\(Y^n_{(k+1)/n^2}< Y^n_{k/n^2}|\cF_k\)
\\
&&\;\;\;\mbox{}\times
\EE\left(\frac{1}{n^2}e^{-\frac{2}{\sqrt{n}}\sum_{l'=1}^{nY^n_{k/n^2}-1}\xi_{l'}\(\WW^{n}_{n^{-2}k}\(Y^n_{k/n^2}\)\)}
\bigg|Y^n_{(k+1)/n^2} <  Y^n_{k/n^2}\,,\  \cF_k\right)\\
&&+\PP\(Y^n_{(k+1)/n^2}> Y^n_{k/n^2}|\cF_k\)\\
&&\;\;\;\mbox{}\times
\EE\left(\frac{1}{n^2}e^{-\frac{2}{\sqrt{n}}\sum_{l'=1}^{nY^n_{k/n^2}}\xi_{l'}\(\WW^{n}_{n^{-2}(k+1)}\(
 Y^n_{k/n^2}+1/n\)\)}
\bigg|Y^n_{(k+1)/n^2}> Y^n_{k/n^2}\,,\ \cF_k\right)\\
&=&\left[\frac{1}{2}-\frac{1}{4\sqrt{n}}\xi_{Y^n_{k/n^2}}\left(
\hat\WW^n_{n^{-2}k}\right) \right]
\frac{1}{n^2}e^{-\frac{2}{\sqrt{n}}\sum_{l'=1}^{nY^n_{k/n^2}-1}\xi_{l'}(\hat{\WW}^{n}_{n^{-2}k})}\\
&&\mbox{}+\left[\frac{1}{2}+\frac{1}{4\sqrt{n}}\xi_{Y^n_{k/n^2}}\left(
\hat\WW^n_{n^{-2}k}\right) \right]
\EE_{\eta}\(\frac{1}{n^2}e^{-\frac{2}{\sqrt{n}}\sum_{l'=1}^{nY^n_{k/n^2}}\xi_{l'}(\hat{\WW}^{n}_{n^{-2}k}+\eta_{1/n})}\)\\
&&\mbox{}+\frac{1}{n^2}1_{Y^n_{k/n^2}=0}
+1_{Y^n_{k/n^2}=K_1}\frac{1}{n^2}e^{-\frac{2}{\sqrt{n}}\sum_{l'=1}^{nK_1-1}\xi_{l'}(\hat{\WW}^{n}_{n^{-2}k})}\,.
\end{eqnarray*}
Therefore,
\begin{eqnarray*}
\lefteqn{\EE(V_{k+1}^2|\cF_k)}\\
&=&\frac{1}{2n^2}\left[\EE_{\eta}\(e^{-\frac{2}{\sqrt{n}}\sum_{l'=1}^{nY^n_{k/n^2}}\xi_{l'}(\hat{\WW}^{n}_{n^{-2}k}+\eta_{1/n})}\)
+ e^{-\frac{2}{\sqrt{n}}\sum_{l'=1}^{nY^n_{k/n^2}-1}\xi_{l'}(\hat{\WW}^{n}_{n^{-2}k})}\right]\\
&&\mbox{}+\frac{1}{4n^{5/2}}\xi_{Y^n_{k/n^2}}\left(
\hat\WW^n_{n^{-2}k}\right)
\left[-e^{-\frac{2}{\sqrt{n}}\sum_{l'=1}^{nY^n_{k/n^2}-1}\xi_{l'}(\hat{\WW}^{n}_{n^{-2}k})}
\right.\\
&&\left.\;\;\;\;\; \mbox{}+
\EE_{\eta}\(e^{-\frac{2}{\sqrt{n}}\sum_{l'=1}^{nY^n_{k/n^2}}\xi_{l'}(\hat{\WW}^{n}_{n^{-2}k}+\eta_{1/n})}\)
\right]
\\
&&\mbox{}+\frac{1}{n^2}\cdot {\bf 1}_{\{Y^n_{k/n^2}=0\}}
+\frac{1}{n^2}\cdot {\bf 1}_{\{Y^n_{k/n^2}=K_1\}}
e^{-\frac{2}{\sqrt{n}}\sum_{l'=1}^{nK_1-1}\xi_{l'}(\hat{\WW}^{n}_{n^{-2}k})}
\\
&=&\frac{1}{2n^2}\left(\EE_{\eta}\(e^{-2B_{Y^n_{k/n^2}}(\hat{\WW}^{n}_{n^{-2}k}+\eta_{1/n})}\)+
e^{-2B_{Y^n_{k/n^2}-1/n}(\hat{\WW}^{n}_{n^{-2}k})}\right)\\
&&\mbox{}+\frac{1}{4n^{5/2}}\xi_{Y^n_{k/n^2}}\left(
\hat{\WW}_{n^{-2}k}\right) \left(-e^{-2B_{Y^n_{k/n^2}-1/n}(\hat{\WW}^{n}_{n^{-2}k})}+
\EE_{\eta}\(e^{-2B_{Y^n_{k/n^2}}(\hat{\WW}^{n}_{n^{-2}k}+\eta_{1/n})}\)
\right)
\\
&&\mbox{}+\frac{1}{n^2}\cdot {\bf 1}_{\{Y^n_{k/n^2}=0\}}
+\frac{1}{n^2}\cdot {\bf 1}_{\{Y^n_{k/n^2}=K_1\}}
e^{-2B_{K_1-1/n}(\hat{\WW}^{n}_{n^{-2}k})}\\
&=&J_{1,n,k}+J_{2,n,k}+J_{3,n,k}+J_{4,n,k}\,.
\end{eqnarray*}
Using the bounds from the proof of Lemma~\ref{lem:5.1} it is easy to see that
\begin{eqnarray*}
 \sum_{k=1}^{\lfloor n^2 t\rfloor}(J_{2,n,k}+J_{3,n,k}+J_{4,n,k}) \rightarrow 0,
\end{eqnarray*}
as $n\rightarrow \infty$. As for $J_{1,n,k}$, again using the convergence
of $(\bW^n, Y^n)$ and the continuity of $B$, it is easy to see that
\begin{eqnarray*}
\sum_{k=1}^{\lfloor n^2t\rfloor} I_{1,n,k} \rightarrow
\int_0^t e^{-2B_{Y_s}(\hat{\WW}_{s})}ds,\;\;{\rm as}\;n\rightarrow \infty,
 \end{eqnarray*}
and we are done.
\qed

\begin{cor}
\label{cor:5.1}
As $n\rightarrow \infty$,
 $M^n$ converges to a continuous 
local
martingale $M$ such that
\begin{eqnarray}
\langle M_{\cdot}\rangle_t &=& \int_0^t e^{-2B_{Y_s}(\hat{\WW}_{s})}ds,\;\;t\geq 0.
\end{eqnarray}
\end{cor}
\paragraph{Proof:}
The continuity of $M$ is immediate from the continuity of the
limiting process $Y$ and Lemma~\ref{lem:5.1}. The rest is immediate
from Lemma~\ref{lem:5.2}.
\qed

\begin{cor}
\label{cor:5.2}
There exists a Brownian motion $\beta$ such that
\begin{eqnarray}
M_{t} &=& \int_0^t e^{-B_{Y_s}(\hat{\WW}_{s})}d\beta_s,\;\;t\geq 0.
\end{eqnarray}
\end{cor}
\paragraph{Proof:} Immediate from the previous corollary.

\paragraph{Proof of Theorem~\ref{thr:5.1}:} Immediate from Lemma~\ref{lem:5.1},
 Corollary~\ref{cor:5.1} and Corollary~\ref{cor:5.2}.
\qed

Finally, we describe the snake process when $g$ is constant. The
description for the general case, more specifically, the uniqueness
of the solution for the martingale problem (\ref{eq:5.1}) remains a
challenging {\em open} problem.

When $g$ is constant, say $g=1$, we have
that
$B_t(x)=B_t$ is a Brownian
motion with constant drift $\nu$. It follows from the martingale
problem (\ref{eq:5.1}) that
\[\int^{Y_t}_0e^{-B_r}dr=\ell^0_t-e^{-B_{K_1}}\ell^{K_1}_t+\int^t_0e^{-B_{Y_s}}d\be_s.\]
Therefore, $Y_t$ is the Brox diffusion reflected at $0$ and $K_1$
(see the Appendix for a description when $\nu=0$).

Next, we consider the conditional (given the lifetime process) path
process. 
Let $w=({\mathbf w},\zeta_w)$ be an element in $\cW$. 
Fix $a\in[0,\zeta_w]$ and $b\ge a$. 
 Similar to LeGall (\cite{LeG99}, p54), we
define $R_{a,b}(w,dw')$ as the unique probability measure on $\cW$
such
that\\
(i) $\zeta_{w'}=b,\;\;R_{a,b}(w,dw')$ a.s.\\
(ii) $w'(t)=w(t)$ for all $t\le a$, $R_{a,b}(w,dw')$ a.s.\\
(iii) Under $R_{a,b}(w,dw')$, $({\mathbf w}'(a+t):\;t\in[0,b-a])$ is a
Brownian motion.

Denote the time set $Q_n=\{n^{-2}k:\;k=0,1,2,\cdots\}$. From the
construction of the discrete snake, it follows that
$\WW^n_s,\;\;s\in Q_n$ is a conditional (given $Y^n$) Markov chain
with transition probability
\[R_{m^n(s,s'),Y^n(s')}(w,dw'),\qquad s<s'\in Q_n,\]
where $m^n(s,s')=\inf\{Y^n(r):\;r\in[s,s']\cap Q_n\}$.

Taking $n\to\infty$, we see that the limit $\{\WW_s,\;s\ge 0\}$ is a
conditional (given $Y$) Markov process with transition probability
\[R_{m(s,s'),Y(s')}(w,dw'),\qquad s<s',\]
where $m(s,s')=\inf\{Y(r):\;r\in[s,s']\}$. Namely, it has the same
conditional law as LeGall's Brownian snake.

\section{Appendix: Convergence to a reflected Brox diffusion}

\setcounter{equation}{0}

We provide in this appendix a short, direct proof of Corollary
\ref{cor-brox} that bypasses the study of the branching process,
relying instead on an embedding of
a random walk in
random environment (RWRE)  into a diffusion in random environment,
in the spirit of \cite{Shi}\footnote{While revising this paper we learnt from
F. Comets about  the paper
\cite{Sei}, that contains a very similar argument for
convergence to the Brox diffusion.}. For backround on
Brownian motion in random environments we refer to
\cite{B},  \cite{Sch}, \cite{T} and to the nice overview in
\cite{Shi}. Background for RWRE can be found in \cite{Zei}.

Recall that a
 Brownian motion in random
environment (BMRE)
is  a process $X_t$ given by
\begin{equation}
\label{eq-oo11} dX_t=d\be_t-\frac12 V'(X_t)dt,
\end{equation}
where $\be_t$ is a Brownian motion and $V$ is called the random
potential. When $V$ is itself a Brownian motion independent of
$\be$, this (formal) process is the Brox diffusion \cite{B}.

We need to
consider reflecting BMRE's.
Let $h$ be the periodic function with period $2K_1$ and $h(x)=|x|$
for $|x|\le K_1$. Let $V$ be a Brownian motion on $x\in[0,K_1]$
and set $\hat V(x)=V(h(x))$ for $x\in\RR$. 
Set formally
\begin{equation}
\label{eq-oo1} dZ_t=d\be_t-\frac12 \hat V'(Z_t)dt.
\end{equation}
(In case $V$ is not smooth, a precise meaning is given to
(\ref{eq-oo1}) by the procedure described in \cite[Section 2]{Shi}).
Let $Y_t=h(Z_t)$.  A {\it formal} application of the
It\^o-Tanaka formula
yields
\begin{eqnarray}
  \label{isitfinal}
dY_t&=&h'(Z_t)dZ_t+d\ell^{Y,0}_t-d\ell^{Y,K_1}_t\\
&=&h'(Z_t)d\be_t-h'(Z_t)
\frac12 \hat V'(Z_t)dt
\nonumber
+d\ell^{Y,0}_t-d\ell^{Y,K_1}_t\\
&=&d\tilde{\be}_t-\frac12 V'(Y_t)dt
+d\ell^{Y,0}_t-d\ell^{Y,K_1}_t,
\nonumber
\end{eqnarray}
where $\tilde{\be}$ is a Brownian motion. 
To justify (\ref{isitfinal}), one argues as follows. First, an application
of Ito's formula for Dirichlet processes, see e.g. \cite{follmer},
gives that for any $g$ which is twice differentiable,
and with $Y^g_t=g(Z_t)$,
\begin{equation}
\label{1601}
 d Y^g_t=g'(Z_t)dZ_t+\frac12 g''(Z_t) dt\,.
\end{equation}
Now note that, by definition of the local time as the occupation time density, 
 the local times of $Z$ and $Y$ at levels $0$ and $K_1$ are 
equal up to multiplicative constant $2$. Therefore a  
standard approximation of $h$ by smooth functions $g$, together with~(\ref{1601}), yields 
(\ref{isitfinal}), provided that the local time
$\ell_t^{Z,x}$
of $Z_\cdot$ is jointly continuous in $t$ and $x$, the latter
at $x=0$ and $x=K_1$.
However, $\ell_t^{Z,x}$ is a continuous transformation of the local
time of the Brownian motion $\beta_t$ (see e.g.
Equation (10) in
\cite{andr} for an explicit formula which holds for any environment---not 
necessarily for the two sided white noise), 
and thus is jointly continuous in its
arguments. This yields (\ref{isitfinal}).
Therefore,  $Y_{\cdot}$ is
a reflecting (at $0$ and $K_1$) Brox diffusion.

\subsection{Embedding}

In this subsection, we introduce an environment and represent $Y^n$
as a RWRE, which we then proceed (after scaling of the environment)
to embed in a diffusion in random environment.

Let the environment be given by a family $\{\xi^n(i),\;i\in\ZZ_+\}$ of
independent random variables with mean $0$
and variance 1. We further
assume that $|\xi^n(i)|\le\sqrt{n}$. 
Define the potential $V^n(\cdot)$ on $\RR_+$ by
\[V^n(x)=\sum^{[x]}_{i=1}\log\frac{\frac{1}{2}-\frac{1}{4\sqrt{n}}\xi^n(i)}
{\frac{1}{2}+\frac{1}{4\sqrt{n}}\xi^n(i)},\] and 
set $\hat V^n(x)= V^n(nh(x/n))$ and let $\hat Z^n$ be the BMRE with
potential $\hat V^n$. Set 
$Z^n(t)=n^{-1}
\hat Z^n(n^2 t)$.
Define the stopping times ${\sigma}^n_0=0$ and
\[{\sigma}^n_{m+1}=\inf\left\{t>{\sigma}^n_m:\;\;\left|{Z}^n(t)
-{Z}^n({{\sigma}^n_m})\right|=1/n\right\}.\]

By Schumacher's theorem (cf. Schumacher \cite{Sch} and Shi
\cite{Shi}), we have
\begin{lemma}
\label{lem:12_01_2}
Let $\tilde{Z}^n_m=n {Z}^n({{\sigma}^n_m}),\;\;m=0,1,2,\cdots$.
Then
$\tilde{Z}^n$ is a RWRE with
\[\PP^{\xi}\(\tilde{Z}^n_{m+1}=i\pm 1\Big| \tilde{Z}^n_m=i\)
=\frac{1}{2}\pm\frac{1}{4\sqrt{n}}\xi^n(nh(i/n)),\]
where $\PP^{\xi}$ is the probability measure conditioned on the environment $\xi$.
\end{lemma}
The next proposition is crucial for the proof of  Corollary~\ref{cor-brox}. 
\begin{prop}
\label{prop:12_01}
The sequence of processes $\left\{\frac{1}{n}\tilde Z^n_{\lfloor tn^2\rfloor}\,, t\geq 0\right\}_{n\geq 1}$ converges weakly in 
$D_{\RR}[0,\infty)$ to the process $Z$ which satisfies (\ref{eq-oo1}). 
\end{prop}
\begin{rem}
\label{rem:12_01}
Note that
$\tilde Y^n\equiv h(\tilde{Z}^n)$
is a sequence of
reflecting (at $0$ and $nK_1$)  RWRE such that  
\[\PP^{\xi}\(\tilde{Y}^n_{m+1}=i\pm 1\Big| \tilde{Y}^n_m=i\)
=\frac{1}{2}\pm\frac{1}{4\sqrt{n}}\xi^n(i),\;\; i=1,\ldots, nK_1-1,\]
and hence by the continuity of the function $h$ and the discussion in the beginning of the appendix, in order to prove Corollary~\ref{cor-brox} it is sufficient to prove Proposition~\ref{prop:12_01}. 
\end{rem}
The rest of the appendix is devoted to the proof of Proposition~\ref{prop:12_01}. 

The following is a
straight-forward consequence of Section 3 of \cite{Shi}.
\begin{lemma}
$Z^n$ is the Brownian motion in random environment with potential
$\hat V^n(nx)$.
\end{lemma}
\paragraph{Proof:}
Let
\[\hat{A}^n_x=\int^x_0e^{\hat V^n(y)}dy.\]
As $\hat {Z}^n$ is the BMRE with potential $\hat V^n$, it is well-known
(see (2.3) in \cite{Shi}) that $\hat{A}^n_{\hat Z^n(t)}$ is a
local martingale with quadratic variation $\hat {\Th}^n(t)$ such
that
\[(\hat {\Th}^n)^{-1}(t)=\int^t_0e^{-2\hat V^n(\hat A^n_{\hat Z^n(u)})}du.\]
We now rescale. Let \[A^n_x=\int^x_0e^{\hat V^n(ny)}dy.\] Then
\[\hat A^n_{\hat Z^n(t)}=\int^{nZ^n\(n^{-2}t\)}_0e^{\hat V^n(z)}dz
=n\int^{Z^n\(n^{-2}t\)}_0e^{\hat V^n(ny)}dy=nA^n_{Z^n(n^{-2}t)}.\] Thus
$A^n_{Z^n(t)}$  is a local martingale with quadratic variation
process $\Th^n(t)=n^2\hat \Th^n(n^2t)$. Thus,
$$
(\Th^n)^{-1}(t)=n^{-2}\int^{n^{-2}t}_0e^{-2\hat 
V^n(nA^n_{Z^n(n^{-2}u)})}du
=\int^t_0 e^{-2\hat V^n(nA^n_{Z^n(u)})}du.
$$
Therefore (see again (2.3) and (2.5) in \cite{Shi}), $Z^n$ is the
BMRE with potential $\hat V^n(nx)$.
\qed

\subsection{Scaling limit}
As was proved in the previous subsection
(see Lemma~\ref{lem:12_01_2}), the scaled RWRE is related
to BMRE by
\begin{eqnarray}
\label{eq:12_01_2}
\frac{1}{n}\tilde{Z}^n_{[n^2t]}=Z^n({{\sigma}^n_{[n^2t]}}).
\end{eqnarray}
 In this
section, we first prove that
\begin{eqnarray}
\label{eq:12_01_3}
{\sigma}^n_{[n^2t]}\to t,\;\; {\rm as}\; n\rightarrow\infty, 
\end{eqnarray} by the strong law of large numbers.  Then, we prove that the
scaled potential for $Z^n$ converges to
$\hat V$, and hence  
$Z^n$ converges to a  
BMRE with potential $\hat V$. This by~(\ref{eq:12_01_2}) and (\ref{eq:12_01_3}) will provide the proof of Proposition~\ref{prop:12_01}.

\begin{lemma}
\label{lem:12_01}
As $n\to\infty$, we have
\[{\sigma}^n_{[n^2t]}\to
t,\qquad a.s.\] uniformly
on compact sets.
\end{lemma}
\paragraph{Proof:} By Proposition 3.2 in \cite{Shi} (or a direct computation
involving a time change), we see that 
$\theta_i=n^2({\sigma}^n_i-{\sigma}^n_{i-1})$, $i=1,2,\cdots$,
are i.i.d. with the same distribution as
\[\theta=\inf\{t>0:\;|W(t)|=1\},\]
where $W$ is a standard Brownian motion. Note that $\EE\theta=1$. By
the strong law of large numbers, we get that
\[{\sigma}^n_{[n^2t]}=t\frac{1}{n^2t}\sum^{[n^2t]}_{i=1}
\theta_i\to t\] uniformly on compacts.
 \qed

For
the next lemma, recall that $Z$ is the processes that satisfies~(\ref{eq-oo1}). 
\begin{lemma}
\label{lem-broxmain} As $n\to\infty$, $Z^n\Longrightarrow Z$ weakly
in $C_{\RR}[0,\infty)$.
\end{lemma}
\paragraph{Proof:}
First we consider the weak convergence of $V^n(nx)$.
Note that
\[V^n(nx)=\sum^{[nx]}_{i=1}\frac{1}{\sqrt{n}}\xi^n(i)+o(1)\equiv
M^n_x
+o(1).\]
Regarding $x$ as the time-parameter, $\{M^n_x,\; x>0\}$ is a
martingale with predictable quadratic variation process
\[\<M^n\>_x=\sum^{[nx]}_{i=1}\EE\(\frac{1}{\sqrt{n}}\xi^n(i)\)^2\to
x\] uniformly on the compacts. Thus, by Theorem 4.13 (\cite{JS},
P358), $M^n$ converges weakly in
$D_{\RR}[0,\infty)$ to a Brownian motion ${V(x), x\geq 0}$.
By
switching to another probability space if necessary, we may 
and will assume that all weak convergences hold almost surely. Then we get, 
$V^n(n\cdot)\to V$, a.s..
Note that by the continuity of $h$, we immediately get that 
\[ \hat{V}^n(x) \rightarrow \hat V(x)= V(h(x)),\;\; {a.s..},\]
uniformly on the compacts of $\RR_+$. 
Note that (see (2.6) in \cite{Shi}),
\[Z^n(t)=(A^n)^{-1}({W^n((T^n)^{-1}(t))}),\]
where
\[A^n_x=\int^x_0e^{\hat V^n(ny)}dy,\]
\begin{eqnarray*}
T^n(t)&=&\int^t_0e^{-2\hat V^n\(n(A^n)^{-1}_{W^n(u)}\)}du,
\end{eqnarray*}
and $W^n$ is a Brownian motion. Since $W^n$ trivially converges 
weakly to the Brownian motion $W$, we assume as before that the
 convergence holds a.s.. 
Then we have
\[A^n_x\to\int^x_0e^{\hat{V}(y)}dy=A_x,\qquad\mbox{ as }n\to\infty, \;\;\]
and
\begin{eqnarray*}
T^n(t)
&\to&\int^t_0e^{-2\hat{V}\(A^{-1}_{W(u)}\)}du=T(t),\qquad\mbox{as
}n\to\infty.
\end{eqnarray*}
Note that all the convergence above are a.s. and uniform on
compacts. We see that
\begin{equation}\label{conv4b} Z^n(t)\to
A^{-1}_{W(T^{-1}(t))}\equiv Z(t).\end{equation}

 By stochastic calculus as in Section 2 of \cite{Shi},
it follows that (\ref{conv4b})
 defines a BMRE $Z(t)$ with potential $\hat V$. \qed

Now
Proposition~\ref{prop:12_01} follows from Lemmas~\ref{lem:12_01}, \ref{lem-broxmain}, and~(\ref{eq:12_01_2}).
Then as we have mentioned already in Remark~\ref{rem:12_01}, Corollary~\ref{cor-brox} follows immediately from Proposition~\ref{prop:12_01}.


\begin{thebibliography}{999}

  \bibitem{andr} P. Andreoletti and R. Diel (2010),
    Limit law of the local time for Brox' diffusion.
    {\em J. Theor. Probab.}, online first, 
    DOI: 10.1007/s10959-010-0314-7.

 \bibitem{B}
T. Brox (1986).
 A one-dimensional diffusion process in a Wiener medium.
{\em Ann. Probab.  \bf 14},  no. 4, 1206--1218.

\bibitem{bib:daw91}
D.~Dawson (1993).
\newblock Measure-valued {P}rocesses.
\newblock {\em \'{E}cole d'\'{e}t\'{e} de Probabilit\'{e}s de Saint Flour},
1991. Lecture notes in Mathematics, {\bf 1541}, Springer, Berlin.


\bibitem{D08}
J.-F. Delmas (2008).
Height process for super-critical continuous state branching
              process.
          {\em Markov Process. Related Fields \bf 14},  309--326.

\bibitem{DS}
 J. Dhersin and L. Serlet (2000).
A stochastic calculus approach for the Brownian snake. {\em Canad.
J. Math.  \bf 52},  no. 1, 92--118.

\bibitem{follmer} H. F\"{o}llmer (1981), Calcul d'It\'{o} sans probabilit\'{e}s.
  {\em S\'{e}m. Prob. XV.},  Lecture Notes in Mathematics  {\bf 850}.
  Springer, Berlin, 143--150.
  
\bibitem{JS}
J. Jacod and A.N. Shiryaev (1987). {\em Limit Theorems for
Stochastic Processes}. Springer, Berlin.

\bibitem{LeG96}
J.-F. Le Gall (1996).
{\em Superprocesses, Brownian snakes and partial differential equations}.
Lecture Notes from the 11th winter school on Stochastic processes,
Sigmundsburg,
Pr\'{e}publication 337 du Laboratoire de Probabilit\'{e}s, Universit\'{e}
Paris VI.


\bibitem{LeG99}
J.-F. Le Gall (1999).
{\em Spatial branching processes, random snakes and partial
              differential equations}.
Lectures in Mathematics ETH Z\"urich, Birkh\"auser Verlag, Basel.





  \bibitem{M}
  L. Mytnik (1996).
  Superprocesses in random environments.
  {\em Ann. Probab. \bf 24}, No. 4. 1953--1978.

\bibitem{Perkins}
E. Perkins (2002).
Dawson--Watanabe superprocesses and measure-valued diffusions.
{\em Lectures on Probability Theory and Statistics, Saint-Flour 1999},
Lecture notes in Mathematics,  {\bf 1781}, Springer, Berlin,
132--329.

\bibitem{Sch}
S. Schumacher (1985). Diffusions with random coefficients. {\em
Contemp. Math. \bf 41}, 351--356.

\bibitem{Sei}
P. Seignourel (2000). Discrete schemes for processes in random media.
{\em Prob. Th. Rel. Fields} {\bf 118}, 
293--322.

\bibitem{Shi}
Z. Shi (2001). Sinai's walk via stochastic calculus. In 
Milieux al\'{e}atoire, {\em Panorama et Synth\'{e}se \bf{12}}, Soc. Math.
France, 53--74.

\bibitem{Si}
Ya.G. Sinai (1982). The limit behavior of a one-dimensional random
walk in a random environment. {\em Th. Probab. Appl. \bf 27},
256--268.

 \bibitem{T}
 H. Tanaka (1995).
 Diffusion processes in random environments.
 {\em Proceedings of the International Congress
 of Mathematicians}, Birkh\"auser, 1047--1054.


\bibitem{Zei}
O. Zeitouni (2004). Random walks in random environment. {\em
Lectures on
probability theory and statistics}, { Lecture Notes in
Mathematics, {\bf 1837}}, Springer, Berlin, 189--312.



\end{thebibliography}
\end{document}